\documentclass[12pt,twoside]{article}
\usepackage{amssymb}
\usepackage[mathscr]{eucal}
\usepackage{eufrak}
\usepackage{amsmath,amsthm}
\usepackage{mathrsfs}
\usepackage[active]{srcltx}
\usepackage[colorlinks=true,
   urlcolor=blue,           filecolor=green,      
   citecolor=green,      
   linkcolor=red,           bookmarks=true,
  unicode,   plainpages=false,   ]{hyperref}
\usepackage{color}
\newcommand{\mg}{\color{magenta}}
\renewcommand{\mg}{\color{black}}

\def\qs{\lfloor\hskip-0.16cm\lceil}
\def\qd{\rceil\hskip-0.16cm\rfloor}

\def\n{\nabla}

\def\intl#1{\int\limits_{#1}}
\def\intll#1#2{\int\limits_{#1}^{#2}}
\def\dm{|\hskip-0.05cm|}
\def\OO{\Omega}
\def\displ{\displaystyle}
\def\VSE{\vspace{6pt}\\&\displ }
\def\VS{\vspace{6pt}\\\displ }
\def\rf#1{{\rm(\ref{#1})}}

\def\R{\Bbb R}
\def\N{\Bbb N}
\def\vep{\varepsilon}
\def\be{\begin{equation}}
\def\ba{\begin{array}}
\def\ea{\end{array}}
\def\ee{\end{equation}}
\def\ov{\overline}

\def\po{\partial\Omega}
\def\bp{\begin{proof}}
\def\ep{\end{proof}}
\font\sc=cmcsc10
\title{\bf\large Global existence of solutions to 2-D Navier-Stokes flow with non-decaying initial data in half-plane\thanks{The research of P.M. was partially supported by GNFM (INdAM) and by MIUR via the PRIN 2017 {\it ``Hyperbolic Systems of Conservation Laws and Fluid Dynamics: Analysis and Applications ''}.
The research of S.S. was partially supported by
JSPS Grant-in-Aid for Scientific Research (B) - 16H03945, MEXT. The latter grant supported a visit of P.M. in Kyoto University. {\mg The authors declare no conflicts of interest.} }}
\author{\sc   Paolo Maremonti\thanks{
Dipartimento di Matematica e Fisica, Seconda
Universit\`{a} degli Studi di
 Napoli, via Vivaldi 43, 81100 Caserta,
 Italy. {\it E-mail address}: 
 paolo.maremonti@unicampania.it} \; and \;Senjo Shimizu
\thanks{{\mg Corresponding author.} Graduate School of Human and Environmental Studies, 
Kyoto University 
Yoshida-nihonmatsu-cho, Sakyo-ku, 
Kyoto 606-8501, Japan. {\it E-mail address}:
{shimizu.senjo.5s@kyoto-u.ac.jp}}}
\begin{document}\maketitle
\newtheorem{ass}
{\bf Assumption}
\newtheorem{defi}
{\bf Definition} 
\newtheorem{tho}
{\bf Theorem} 
\newtheorem{rem}
{\sc Remark} 
\newtheorem{lemma}
{\bf Lemma} \newtheorem{coro}
{\bf Corollary} 
\newtheorem{prop}
{\bf Proposition}
\renewcommand{\theequation}{\arabic{equation}}
\setcounter{section}{0}\par {\small {\bf Abstract} - We investigate  the Navier-Stokes initial boundary value problem in the half-plane  $\R^2_+$ with   initial data $u_0\in L^\infty(\R^2_+)\cap J^2_0 (\R^2_+)$ or with non decaying initial data $u_0\in L^\infty(\R^2_+)\cap J^p_0 (\R^2_+),\,p>2\,.$ We introduce a technique that allows to solve the two-dimesional problem, further, but not least, it can be also employed to obtain weak solutions, as regards the non decaying initial data, to the  three-dimensional Navier-Stokes IBVP. This last result is the first of its kind.
\par\vspace{0.5pc}
{\mg {\it Keywords}.    2-D Navier-Stokes equations, non-decaying data, global solution}
}
\section{\label{In}Introduction}In this paper   we consider the following Navier-Stokes initial boundary value problem
\be\label{NS}\ba{l}u_t+u\cdot\n u+\n\pi_u=\Delta u\,,\quad \n\cdot u=0\mbox{ in }(0,T)\times\OO\,,\VS u=0\mbox{ on }(0,T)\times\po,\quad u=u_0\mbox{ on }\{0\}\times\OO\,,\ea\ee where the symbol $\OO$   denotes an exterior domain,  $\R^n$ and $\R^n_+$, $n\geq2$, and by $a\cdot\n b$ we mean $(a\cdot\n)b$. We look for solutions global in the time to problem \rf{NS} with non decaying initial data.   The problem  of the existence of solutions to \rf{NS} with non decaying data   has been considered by several authors and, we think that the first results, where $n\geq2$,  go back to the papers \cite{GMZ, GIM, GMS, MR1,  MR2}. But the special case of the two-dimensional problem involves a particular interest for the possibility to obtain global existence in the pointwise norm. A natural setting of the problem is the function space $L^\infty((0,T)\times\OO)$. In this sense a first result is given by Giga, Matsui and Sawada in \cite{GMS} limited to the Cauchy problem. Subsequently, in \cite{ST} Sawada and Taniuchi improve the $L^\infty$-norm of the solutions of \cite{GMS}. Based on a result by Zelik in \cite{Z}, a recent contribute given by Gallay in \cite{Gl} establishes an estimate that up today is the best one:
$$\dm u(t)\dm_{L^\infty(\R^2)} \leq c\dm u_0\dm_{L^\infty(\R^2)}(1 + c\dm u_0\dm_{L^\infty(\R^2)}t), \mbox{ for all }t > 0.$$ However all these results concern the Cauchy problem associated to the 2D-Navier-Stokes equations with non decaying data. Subsequently, the problem has been considered in exterior domains. Firstly Abe in \cite{ABE} gives a result of local existence of the mild solution with initial data $u_0\in L^\infty(\OO)$  that can be seen as a weak solution to the Navier-Stokes problem. Then, in \cite{MSh} Maremonti and Shimizu   improve the result by Abe giving the existence and uniqueness of solutions  to the   Navier-Stokes initial boundary value problem in exterior domains which are defined for all $t>0$. Actually, these authors are able to prove a smooth extension of the solution determined by Abe. The  results contained  in \cite{MSh} can be also seen as a ``structure theorem'' of the weak solution given in \cite{ABE}. The result by Maremonti and Shimizu is based on the possibility to reduce the problem to an $L^2$-theory. In the sense that the solution $u$ is seen as the sum of three fields, that is $u=U+W+w$, where $U$ and $W$ are solutions to a linear problem and  keep the non decaying character of the initial data, instead $w$ is the  solution to a nonlinear perturbed Navier-Stokes with homogenous initial data and suitable  force data  with compact support. For the field $w$ is applicable the $L^2$-theory (see e.g. \cite{La}). However this approach seems to be unable to work in the case 
 of $\po$ not bounded. \par More recently,   in \cite{ABED}, as particular case of the results by Maremonti and Shimizu,  Abe proves global existence in exterior domains by means of  the special assumption of $u_0\in L^\infty(\OO)$ $(\OO\subseteq\R^2$) and $\dm \n u_0\dm_2<\infty$. In the case of the half-plane, he obtains a result under the assumption that the initial data is decaying 
from the viewpoint of Hardy's inequality.  \par  Although the geometry of the half-plane, and more in general the one  of the half-space, concerns a particular case of the mathematical theory, it is very interesting in the applications and recall the attention of several authors \cite{CJ-III, CJ-I, CJ-II, HM}. Therefore the aim of the present paper is to prove that the result obtained by Maremonti and Shimizu in \cite{MSh} also holds in the half-plane. In order to state our chief results we introduce some notations. \par 
 By the symbol $\mathscr C_0(\OO)$, we denote the set of all solenoidal vector fields $\varphi\in C_0^\infty(\OO)$.
 By the symbol $J^q(\OO)$, $q\in(1,\infty)$, we indicate  the completion of $\mathscr C_0(\OO)$ in $L^q(\OO)$ Lebesgue space. The symbol $P_q$ indicates the projector from $L^q(\R^2_+)$ into $J^q(\R^2_+)$. We set $J^{1,q}(\OO)$:=completion of $\mathscr C_0(\OO)$ in Sobolev space $W^{1,q}(\OO)$. We set  $J^{q}_0(\R^2_+):=$ completion of $\mathscr C_0(\R^2_+)$ in the  space $\dot W^{1,q}(\R^2_+)$, that is with respect to the seminorm $\dm \n\cdot\dm_q\,.$
\begin{tho}\label{CT}{\sl Let $u_0\in L^\infty(\R^2_+)\cap J_{0}^p(\R^2_+)$, $p\in(2,\infty)$. Then there exists a unique solution $(u,\pi_u)$ to problem \rf{NS} such that \be\label{CT-I}\ba{ll}\mbox{for all }\eta>0\mbox{ and }T>0,\hskip-0.2cm& u\!\in\! C([0,T)\times\R^2_+)\cap C^2((\eta,T)\times\R^2_+)\,,\VSE u_t,\n\pi_u\!\in\! C((\eta,T)\times \R^2_+)\,.\ea\ee Moreover, up to a function $c(t)$, we get the estimate \be\label{PFE} |\pi(t,x)|\leq P(t)|x|^{\mu},\mbox{ for all }t>0\mbox{ and }x\in \R^2_+\,,\ee for a suitable $\mu\in(0,1)$, and, for all $\vep>0$ and $T>0$, $P(t)\in C(\vep,T)\cap L^s(0,T)$ for a suitable $s>1$.}\end{tho}
\begin{rem}{\rm Comparing the assumptions made on $u_0$   in \cite{MSh} with the one of Theorem\,\ref{CT}, we point out that
other the integrability of $\n u_0$, as a consequence of the Sobolev embedding, we  assume more regularity for the initial data.\par
 Although it is possible  to deduce   an estimate of the $\dm u(t)\dm_\infty$ for all $t>0$, for the sake of the brevity, we do not give it. Like in the paper \cite{MSh}, the problem of the existence and the bound of the uniform norm of the solutions are two different questions. \par Unlike all the quoted results,   Theorem\,\ref{CT} enjoys of a quite original proof, which aquires   a further interest for its application to the three-dimensional case. The proof is a variant of the one exhibits in the paper \cite{MDS,MLW}, where the solutions   are decaying in some sense. More precisely, firstly  we introduce a finite family of solutions, each is the solution of a   Navier-Stokes linearezed problem. The first element of the family is the solution to the Stokes problem with non decaying initial data $u_0$. The number of the solutions depends on   the exponent $p$ of $\dm\n u_0\dm_p$. If $p=2$, then we have just one linear (Stokes) problem, hence the solution $(U,\pi_U)$. Then, in order to solve problem \rf{NS}, we study the solution $(w,\pi_w)$ to the perturbed (nonlinear) Navier-Stokes problem where the coefficients are $U$ and $\n U$, the problem has an homogeneous initial data and body force $U\cdot\n U$, that as a matter of course belongs to $L^2$ (see problem \rf{PNS}). So a $L^2$-theory is applicable to prove the existence of $(w,\pi_w)$. Hence $u:=U+w$ solves \rf{NS}. If $p>2$, we consider the greatest integer floor $k$ of $\log_2 \frac p2$. It is such that $\frac p{2^k\!\!\!\null}>2$ and $\frac p{2^{k+1}\!\!\!\!\!\!\!\!\!\null}\;\;\;\leq 2$. For $h\in\{1,\dots,k-1\}$ we consider $(v^h,\pi_{v^h})$ solution of a corresponding linearezed problem (see \rf{LNS}) where the coefficient are $v^\ell,\n v^\ell$, $\ell=1,\dots,h,$ the initial data is zero and the force data is the $v^{h-1}\!\cdot\!\n v^{h-1}\!\in\! L^\frac p{2^h\!\!\!\null}\,(0,T;L^\frac p{2^h\!\!\!\null}\,(\R^2_+))$.\,Since $\frac p{2^h\!\!\!\null}>\!\frac p{2^{h+1}\!\!\!\!\!\!\!\!\!\null}\hskip0.4cm$     by construction the last step is a field $(w,\pi_w)$ solution to a nonlinear perturbed Navier-Stokes problem with homogeneous initial data and force data $F:=-v^k\cdot \n v^k$ that, recalling the definition of $k$, belongs to $ L^2(0,T;L^2(\R^2_+))$ (see problem \rf{KPNS}). In this final step we can apply the $L^2$-theory which allows to conclude the proof of Theorem\,\ref{CT}.}\end{rem} The approach used in the proof of Theorem\,\ref{CT} also allows us to deduce the following result for the 3D-Navier-Stokes Cauchy problem and IBVP in the half-space:
\begin{tho}\label{3D}{\sl Let $u_0\in L^\infty(\OO)\cap J^p_0 (\OO)$, $p>3$ and $\OO\equiv \R^3$ or $\OO\equiv\R^3_+$.     Then there exists a field $u:(0,\infty)\times\OO\to \R^3$ which is a solution in the distributional sense to problem \rf{NS}. Moreover, set $k$    the greatest integer floor of $\log_2\frac p2$, we get $u:=U+{\overset {k+1}{\underset {\ell=1}\sum}}v^\ell+w$, where for all $\eta>0$ and  $T>0$,  $U,v^\ell\in C([0,T)\times\OO)\cap C^2((\eta,T)\times\OO) $, $U_t,v^\ell_t,\n\pi_U,\n\pi_{v^\ell}\in C((\eta,T)\times\OO)$, and   $w\in L^\infty(0,T;J^2(\OO))\cap L^2(0,T;J^{1,2}(\OO))$. Finally, the solution u is strongly continuous to the initial data, $\displ\lim_{t\to0}\dm u(t)-u_0\dm_\infty=0.$}\end{tho} Apart of an interesting, but a special, result obtained by Sawada in \cite{Sw}, as far we know, the result  of Theorem\,\ref{3D} is the first of its kind. We do not give the proof of Theorem\,\ref{3D}. Formally it is quite analogous to the one of the 2D case. The unique difference concerns the last step. Actually, for the field $w$ we have to employ a Hopf-Leray existence theorem. This last gives to our solution the character of weak solution, and makes the difference with the 2D case, for which the $L^2$-theory allows to obtain regular solutions (see e.g. \cite{La}). In the claims of Theorem\,\ref{3D} it has not to surprise that the initial data is assumed continuously with respect to the uniform norm. This is a consequence of the fact that for all data $u_0\in C(\OO)$ at least on some interval $(0,T_0)$ the solution is regular as proved in \cite{GIM,MR1}.\par The authors would like to conclude the introduction giving special thanks to Professor Yasushi Taniuchi who made our attention to the IBVP for the half-plane problem with non-decaying data.
\section{\label{N-P}Some notations and preliminaries}
 By the symbol $\chi$, we denote a smooth positive cutoff function such that $\chi(\rho)=1$ for $\rho\leq 1$, $\chi(\rho)\in [0,1]$ for $\rho\in[1,2]$ and $\chi(\rho)=0$ for $ \rho\geq2$. For $R>0$, we set $\chi_R(x):=\chi(\frac {|x|}R)$.   For $q\ne2$, we set $\ov l:=\big[2-\frac2q\big]$ and by
$\widetilde B^{2-\frac2q,q}(\OO)\subset L^q(\OO) $, $\OO=\R^2$ or $\R^2_+$,   we mean the set of functions such that   \be\label{BS}  \dm u\dm_{2-\frac2q,q}:= \ll u\gg_q^{2-\frac2q}+\dm u\dm_{W^{\ov l,q} }<\infty\,,\ee where the functional $ \ll\cdot \gg^{2-\frac2q}_q$ is given by $$ \ll u \gg_2^{2-\frac2q}=\Big[ \intl{\OO}|z|^{-2q} \Big[\intl{\OO }|u(x)-u(x+z)|^qdx\Big]dz\Big]^\frac1q\,.$$  For $q=2$ we set $\widetilde B^{1,2}(\OO):=W^{1,2}(\OO)$.
\par By the symbol $C^{k,\lambda}(\OO)$, $k\in\N$ and $\lambda\in(0,1)$, we denote the H\"older's space of functions continuous differentiable with their  derivatives $D^\alpha u$, $|\alpha|\leq k$, and with $D^\alpha u,\,|\alpha|=k$ $\lambda$-H\"older continuous. The norm in $C^{k,\lambda}$ is indicated by $\dm\cdot\dm_{k,\lambda}$ and H\"older's seminorm by $\qs\cdot\qd^{(\lambda)}_\OO$. We use the symbol $\qs\cdot\qd^{(\lambda)}$ when there is no confusion about the domain. \par
Let $q\in[1,\infty)$, let $X$ be a Banach
space with norm $\|\cdot\|_X$. We denote by $L^q(a,b;X)$ the set of all
function $g:(a,b)\to X$ which are measurable and such that the
Lebesgue integral \mbox{\footnotesize $\displ\int_a^b$}$\|
g(\tau)\|^q_Xd\tau= \|g\|_{L^q(a,b;X)}^q<\infty$. As well as, if
$q=\infty$ we denote by $L^\infty(a,b;X)$ the set of all function
$g:(a,b)\to X$ which are
 measurable and such that ${\rm ess\ sup}_{t\in(a,b)}\ \|
g(t)\|_X=\|g\|_{L^\infty(a,b;X)}<\infty$.
Finally, we denote by $C([a,b);X)$ (resp. $C(a,b;X)$) the set of functions which are continuous from $[a,b)$ (resp. $(a,b)$) into $X$ and normed with $\displ\sup_{(a,b)}\dm g(t)\dm_X<\infty$.
\begin{lemma}[Sobolev embedding]\label{SOB}{\sl Let $f\in \dot W^{1,q}(\R^2_+),\, q\in[1,2)$. Then there exist a constant $f_0$ such that
\be\label{SOB-I}\dm f-f_0\dm_{\frac{2q}{2-q}}\leq c\dm \n f\dm_q\,,\ee where $c$ is independent of $f$.\,If $f\!\in\! \dot W_0^{1,q}(\R^2_+)$ then inequality \rf{SOB-I} holds with $f_0=0$.}\end{lemma}\bp See e.g. \cite{Ga} Section\,II.5\,.  \ep
\begin{lemma}[Gagliardo-Nirenberg inequality]{\sl Let $u\in L^q(\R^2_+)$ with $\n^m u\in L^r(\R^2_+)$,  $m=1,2$. Assume that $\gamma_{tr}(u)=0$. Then,  there exists a constant $c$ independent of $u$ such that
\be\label{GN-I}\dm \n^ju\dm_s\leq c\dm \n^m u\dm_r^a\dm u\dm_q^{1-a}\,,\ee   provided that, for $j=0,1$ and $m=1,2$, $1\leq q$, $r\leq\infty    $,   $\frac 1s=\frac j2+a(\frac1r-\frac m2)+(1-a)\frac1q\,$ with $a\in[\frac jm,1]$ for $r\ne2$,     with $a\in[\frac 12,1)$ for $r=2$. Further, if $u\in J^q(\R^2_+)$ and $m=2$, then, for $r\in(1,\infty)$,  we also get
\be\label{GN-II}\dm \n^ju\dm_s\leq c\dm P_r\Delta u\dm_r^a\dm u\dm_q^{1-a}\,.\ee   
 }\end{lemma}\bp Inequality \rf{GN-I} is   the well Gagliardo-Nirenberg inequality. Inequality \rf{GN-II} is related to the solenoidal functions and it is consequence of the fact that $\dm \n^2u\dm_r\leq c\dm P_r\Delta u\dm_r$\,\footnote{We remark that in the case of $m=2$ inequality \rf{GN-I} also can be given by means of $\Delta$ in place of $\n^2$.} for $r\in(1,\infty)$.\ep We recall the following well known version of Gronwall inequality: \begin{lemma}[Gronwall lemma]\label{GW}{\sl Assume that $\varphi(t),\,\psi(t),\,h(t)$ and $k(t)$ are continuous and nonnegative functions on $[0,T]$. Assume that the following integral inequality holds:$$\varphi(t)+\intll0t\psi(\tau)d\tau\leq \intll0th(\tau)\varphi(\tau)d\tau+\intll0tk(\tau)d\tau\,\mbox{\; for all }t\in[0,T]\,.$$ Then we get $$\varphi(t)+\intll0t\psi(\tau)d\tau\leq exp\big[\intll0t h(\tau)d\tau\big]\intll0t k(\tau)d\tau\,\mbox{\; for all }t\in[0,T]\,.$$}\end{lemma}
 \section{\label{StPr}Stokes problem}
In this section we study the following initial boundary value problem:
\be\label{ST}\ba{l} U_t+\n \pi_U=\Delta U+G,\quad \nabla \cdot U=0\;\mbox{ in }(0,T)\times\R^2_+\,,\VS U=0\mbox{ on }(0,T)\times \{x_2=0\}\,,\quad U=u_0\mbox{ on }\{0\}\times\R^2_+\,.\ea\ee 
In order to discuss problem \rf{ST} by means of the Green function  with the special assumption $u_0\in L^\infty(\R^2_+)\cap J^{p}_0(\R^2_+),$ we have to premise some results. We start with the following well known result:
\begin{lemma}\label{HIH}{\sl Assume that $u\in L^p_{\ell oc}(\ov\R^2_+)$ with $\n u\in L^p(\R^n_+),\,p\in[1,\infty)$. Assume that $u(x_1,0)=0$ for all $x_1\in\R$. Then, for all $\widetilde R\geq0$, we get
\be\label{HIH-I}\intl{|x|>\widetilde R}|u(y)|^p|y|^{-p}dy\leq \pi^p\!\!\!\intl{|x|>\widetilde R}|\n u(y)|^pdy\,.\ee}\end{lemma}\bp We reproduce the proof for the sake of the completeness. Introduced a polar coordinate frame $(r,\theta)$, almost everywhere in $r>0$, we get
$$\intll0\pi|u(r,\theta)|^pd\theta\leq \pi^p\intll0\pi|\mbox{$\frac\partial{\partial \theta}$}u(r,\theta)|^pd\theta\,.$$
Recalling that $|\frac{\partial u}{\partial \theta}|\leq r|\n u|$, we obtain 
$$\intll0\pi|u(r,\theta)|^pr^{1-p}d\theta\leq \pi^p\intll0\pi|\mbox{$((\frac\partial{\partial r}u(r,\theta))^2+(\frac1r\frac\partial{\partial \theta}$}u(r,\theta))^2|^\frac p2rd\theta\,.$$
Integrating this last inequality on $(\widetilde R,\infty)$, we deduce the thesis.\ep The following   is also a well known result (see e.g. \cite{MSt} and \cite{ABED}). Again for the sake of the completeness we reproduce the proof given in \cite{MSt}.
\begin{lemma}\label{DEN}{\sl Let $u_0\in L^\infty(\R^2_+)\cap J^{p}_0(\R^2_+),\,p\geq2$. Then there exists a sequence $\{u^n_0\}\subset L^\infty(\R^2_+)\cap J^{p}_0(\R^2_+)$,   such that,  for all $n\in\N$, $u_0^n$ has compact support with \be\label{BN}\ba{l}\dm u_0^n\dm_\infty\leq c\dm u_0\dm_\infty\,,\mbox{ for all }n\in\N\,,\VS\dm\n u_0^n\dm_p\leq \dm \n u_0\dm_p+o(1)\,,\mbox{ for all }n\in\N,\ea \ee and for all $R>0$, the sequence converges to $u_0$ in $L^\infty(B_R\cap \R^2_+)\cap J^{p}_0(\R^2_+) $. }\end{lemma}\bp 
We denote by $\{\chi^n\}$ the sequence of cutoff functions with $\chi^n:=\chi(\frac xn)$, where $\chi(\rho)$ is the cutoff function introduced in section\,2. Hence, for all $n\in\N$, the support is the ball $B_{2n}$ and  $|\n \chi^n(x)|\leq \frac cn$. For all $n\in\N$, we set $\ov u^n:=u_0\chi^n$ and we consider the Bogovski problem \be\label{Bg}\ba{l}\nabla \cdot 
\widetilde u^n=-\n\cdot\ov u^n=-\n\chi^n\cdot u_0 \mbox{\; in\; }(B_{2n}-B_n) \cap \R^2_+\,,\VS \widetilde u^n=0\mbox{\; on \;}\partial(({B_{2n}-B_n}) \cap \R^2_+)\, .\ea\ee
Since $u_0$ is divergence free, in problem \rf{Bg} the compatibility condition    is satisfied. Hence there exists at least a solution $\widetilde u^n$, and since the domain $(B_{2n}-B_n) \cap \R^2_+$ is of homothetic kind,   with a constant $c$ independent of $n$   we obtain  \be\label{Bg-I}\ba{l}\dm \nabla \widetilde u^n\dm_p\leq c\dm u_0\cdot\n\chi^n\dm_p\leq cn^{-1+\frac2p}\dm u_0\dm_\infty\,,\VS \dm \n \widetilde u_n\dm_p\leq c\dm \n u_0\dm_{L^p(|x|\geq n)}\,,\ea\ee where for the latter estimate we have employed Lemma\,\ref{HIH}. By the Poincar\'e inequality and \rf{Bg-I}$_1$\be\label{Bg-II}\dm \widetilde u^n\dm_p\leq cn\dm \n \widetilde u^n\dm_p\leq cn^\frac2p\dm u_0\dm_\infty\,.\ee Employing the Gagliardo-Nirenberg inequality, via estimates \rf{Bg-I}$_1$-\rf{Bg-II} we get
\be\label{Bg-III}\dm \widetilde u^n\dm _\infty\leq c\dm u_0\dm_\infty\,\mbox{\; for all }n\in\N.\ee We extend $\widetilde u^n$  to zero on $\R^2_+$.  For all $n\in\N$, the extension is  denoted again by $\widetilde u^n$. We  define $u_0^n:=\ov u^n+\widetilde u^n$.  Hence it follows that  the sequence $\{u_0^n\}\subset L^\infty(\R^2_+)\cap J^{p}_0(\R^2_+)$. Trivially we get that, for all $R>0$, the sequence $\{u^n_0\}$ converges in $ L^\infty(B_R\cap\R^2_+)$. Then Lemma\,\ref{HIH} estimate   \rf{Bg-I}$_2$ ensure the convergence of the sequence $\{u^n_0\}$ in $J^{p}_0(\R^2)$. Actually we get $$\dm \n u_0-\n u^n\dm_p\leq\dm (1-\chi^n)\n u_0\dm_p+\dm u_0\n\chi^n\dm_p+\dm \n\widetilde u^n\dm_p\leq c\dm \n u_0\dm_{L^p(|x|>n)}\,.$$
\ep
We represent the solution to problem  \rf{ST} by means of the Green function furnished in \cite{Sl} and in \cite{Slll}, see also  \cite{U}. In two-dimensional case  the Green function is defined as follow:
\be\label{Gvl}\ba{l}\displ \frak G_{11}(t,x,y)\!:=
  \Gamma(t,x\!-\!y)-\!\Gamma(t,x\!-\!y^*)
 +4  D_{x_1}\!\!\intll0{x_2}\!\!\intl{\R}\!\!D_{x_1}E(x\!-\!z)\Gamma(t,z\!-\!y^*)dz\VS 
\frak G_{12}(t,x,y)\!:=0\VS \frak G_{21}(t,x,y)\!:=
 4  D_{x_1}\!\!\intll0{x_2}\!\!\intl{\R}\!\!D_{x_2}E(x-z)\Gamma(t,z-y^*)dz  \VS \frak G_{22}(t,x,y)\!:=
  \Gamma(t,x-y)-\Gamma(t,x-y^*)
  \ea\ee
\be\label{Gpr}\ba{l} \displ  \frak P_1(t,x,y):=4 D_{x_1}\intl{\R}\Big[D_{x_2}E(x_1-z_1,x_2)\Gamma(t, z_1-y_1,y_2) \VS\hskip 5cm+E(x_1-z_1,x_2)D_{y_2}\Gamma(t,z_1-y_1,y_2)\Big]dz_1\VS \frak P_2(t,x,y):=0\,.\ea  \ee
In formula \rf{Gvl}-\rf{Gpr} function $\Gamma(t,z)$ is the kernel of heat equation and function $E(z)$ is the fundamental solution of Laplace equation  and $y^*:=(y_1,-y_2)$. we denote by $\Gamma^*(t,x,y):=\Gamma(t,x-y^*)$. 
We set $$\frak G^*_{r1}(t,x,y):=4  D_{x_1}\!\!\intll0{x_2}\!\!\intl{\R}\!\!D_{x_r}E(x-z)\Gamma(t,x-y^*)dz\,,\;r=1,2.$$ For all   $k\in\N$,  $h:=(h_1,h_2), \,h_i\in\N,\,i=1,2$, and $l:=(l_1,l_2), \,l_i\in\N,\,i=1,2$, and $\mu>0$,  the following estimates hold:
\be\label{EGF}\ba{l}  |D^{k,h}_{t,z}\Gamma(t,z)|\!\leq ct^{\frac\mu2}(|z|^2\!+t)^{-1 -k-\frac\mu2-\frac{|h|}2},\mbox{ for all }(t,z)\!\in\!(0,T)\times\R^2_+\,,\VS   \Big|D^{k,h,l}_{t,x,y}\frak G^*_{r1}(t,x,y)\Big|\!\leq ct^{-k-\frac {l_2}2 }(x_2^2\!+\!t)^{-\frac{ h_2}2 }{\rm exp}[-\widehat a \frac yt{\hskip-0.1cm\null^{\null^{2\atop 2}}} ](|x-y^*|^2\!+\!t)^{-1-\frac{h_1+l_1}2}   ,\VS |D_{t,x,y}^{k,h,l}\frak P_j(t,x,y)|\leq c t^{-1-k-\frac {l_2}2 } {\rm exp}[-\widehat a \frac yt{\hskip-0.1cm\null^{\null^{2\atop 2}}} ](|x-y^*|^2\!+\!t)^{-\frac{1+h_1+h_2+l_1}2} , \ea \ee with   $c$ and $\widehat a$ positive constants independent of $t,x$ and $y$\,.\par Via suitable hypotheses for data $u_0$, a solution to the Stokes problem \rf{ST} can be represented in components as
\be\label{SRF}\ba{l}\displ U_i(t,x):=\mbox{${\underset{j=1}{\overset 2\sum}}$}\intl{\R^2_+}\frak G_{ij}(t,x,y){u_0}_j(y)dy=:\frak G_i[u_0]\,,\;i=1,2\, ,\VS\pi_U(t,x)\!:=\mbox{${\underset{j=1}{\overset 2\sum}}$}\intl{\R^2_+} \frak P_j(t,x,y){u_0}_j(y)dy=:\frak P[u_0]\,.\ea\ee
\begin{tho}\label{ESP}{\sl Let $G=0$ in \rf{ST}. For all $u_0\in L^\infty(\R^2_+)$ with null divergence, there exists a unique smooth solution to problem \rf{ST} such that
\be\label{ESP-I}\ba{l}\dm U(t)\dm_\infty\leq c\dm u_0\dm_\infty\,,\mbox{ for all }t>0\,,\VS \dm \n U(t)\dm_\infty\leq ct^{-\frac12}\dm u_0\dm_\infty\,,\,\mbox{ for all }t>0\,,\ea\ee with $c$ independent of $u_0$. Moreover,  we get
\be\label{IC}\ba{l}\displ\mbox{if }u_0\in L^\infty(\R^2_+)\,,\mbox{\, then }\lim_{t\to0}\,(U(t)-u_0,\varphi)=0\,\mbox{ for all }\varphi\in \mathscr C_0(\R^2_+)\,,\VS \mbox{if }u_0\in C(\R^2_+)\,,\mbox{\, then }\lim_{t\to0}\dm U(t)-u_0\dm_\infty=0\,.\ea\ee} Finally, up to a function $c(t)$, for the pressure term $\pi_U$ we get
\be\label{PF}\gamma\in(0,1),\quad|\pi_U(t,x)|\leq c\dm u_0\dm_\infty|x|^{\gamma}t^{-\frac\gamma2}\,,\mbox{ for all }(t,x)\in (0,T)\times\R^2_+\,.\ee \end{tho}
\bp In the case of $u_0\in C(\R^2_+)$ see Solonnikov \cite{Sll}, Theorem\,1. In the case of $u_0\in L^\infty(\R^2_+)$, the proof proposed by Solonnikov need some minimum modifications. We have to prove  that the solution $U$  given by representation formula \rf{SRF}$_1$  satisfies the initial condition \rf{IC}$_1$. Further, the proof of the uniqueness, proposed by Solonnikov in \cite{Sll}, also works in the case of $u_0\in L^\infty(\R^2_+)$. Actually the proof is based on an argument of duality, hence the property \rf{IC}$_1$ is sufficient in order to prove the reciprocity formula given by  the duality. We consider the sequence $\{u_0^n\}$ of Lemma\,\ref{DEN}. We have that $\dm u_0^n\dm_\infty\leq c\dm u_0\dm_\infty$, and $u_0-u_0^n$ has a support for $|x|>n$, for all $n\in\N$.  Hence, choosing $n$ sufficienty large such that supp\,$\phi\subset B_n$, we get \be\label{IC-I}(U(t)-u_0,\phi)=(U(t)-U^n(t),\phi)+(U^n(t)-u_0^n,\phi)\,,\ee where $U^n(t)$ is the solution of problem \rf{ST} corresponding to the initial data $u_0^n$. Via the representation formula, employing estimates \rf{EGF}$_{1,2}$, for all $x\in$\,supp\,$\phi$, we get $$\ba{ll}|U(t,x)-U^{n+m}(t,x)|\hskip-0.2cm&\displ=|\frak G[u_0-u^{n+m}_0]|\VSE\leq c\dm u_0\dm_\infty\hskip-0.3cm\intl{|y|>n+m}\hskip-0.4cmt^\frac\mu2|y|^{-2-\mu}dy+\intl{|y| >n+m}\hskip-0.5cm\mbox{$\exp[-\hat a\frac{y_2^2}t]$}|y|^{-2}dy\VSE\leq c\dm u_0\dm_\infty\big[t^\mu(n+m)^{-\mu}+t^\frac12(n+m)^{-1}\big]\,,\ea$$ with $c$ independent of $t$ and $n,m$. Hence for the former term of \rf{IC-I} we obtain the estimate
$$|(U(t)-U^{n+m}(t),\phi)|\leq c\dm u_0\dm_\infty\big[t^\mu(n+m)^{-\mu}+t^\frac12(n+m)^{-1}\big]\dm \phi\dm_1,\mbox{ for all }t>0\,.$$
For the latter term of \rf{IC-I}, applying H\"older's inequality, for $q\in(1,\infty)$, we get
$$|(U^n(t)-u_0^n,\phi)|\leq \dm U^n(t)-u_0^n\dm_q\dm \phi\dm_{q'}\,\mbox{ for all }t>0\,.$$ Therefore for $t\to0$ we obtain \rf{IC}$_1$. The estimate \rf{PF} is contained in \cite{MR1}, Theorem\,2.1\,. \ep
\begin{tho}\label{GPST}{For all $u_0\in L^\infty(\R^2_+)\cap J^{p}_0(\R^2_+)$, $p\in(1,\infty)$,    the solution furnished by Theorem\,\ref{ESP} verifies \be\label{GPST-I}\ba{l}\dm \n U(t)\dm_p\leq c\dm \n u_0\dm_p\,,\mbox{ for all }t>0\,,
\ea\ee with a constant $c$ independent of $u_0$.}\end{tho}
\bp Firstly we assume $p\geq2$. By virtue of Lemma\,\ref{DEN} there exist a sequence $\{u_0^n\}\subset L^\infty(\R^2_+)\cap J^{p}_0(\R^2_+)$ converging to $u_0$ in $L^\infty(B_R\cap\R^2_+)\cap J^{p}_0(\R^2_+)$, where $u_0^n:=\ov u^n+\widetilde u^n$ has a compact support, and enjoys estimates \rf{BN}. We denote by $(U^n,\pi_{U^n})$ the sequence of solutions corresponding to $\{u_0^n\}$, where $U^n:=\frak G[u^n_0]$ and $\pi_{U^n}:=\frak P[u^n_0]$.  Recalling that $u_0-u^n_0=(1-\chi^n)u_0+\widetilde u^n$ has a support for $|x|>n$, by representation formula and estimates \rf{EGF}, for all $t>0$, we get that $$\ba{l}\displ|U(t,x)-U^n(t,x)|\leq c\dm u_0\dm_\infty \Big[t^{\frac\mu2}\!\!\!\!\intl{|y|>n}\!\!\!(|x-y|^2+t )^{-1-\frac\mu2 }dy\VS\hskip6cm+c\!\!\!\intl{|y|>n}\!\!\!\! e^{-    \mbox{$\frac{\widehat a}t$}\hskip0.05cm{{y ^2_2}}}(|x-y^*|^2+t )^{-1}dy\Big]\,,\ea$$   
which ensures the uniform convergence on any compact $K$ subset of $\R^2_+$. Analogous is the proof of the following convergence:
\be\label{PGE}\mbox{for all }t>0,\mbox{ and compact }K\!\subset\!\R^2_+\!:\lim_{n\to\infty}\dm \n U(t)\!-\!\n U^n(t)\dm_{L^\infty(K)}\!=0\,.\ee 
From the representation formula of $U^n$ it is not difficult to deduce for the tangential derivative that
$$i=1,2,\quad D_{x_1}U_i^n=\frak G_{ij}[D_{y_1}{u_0^n}_j]\,, \mbox{ for all }(t,x)\mbox{ and }n\in\N\,. $$
From the equation of the divergence we also get $$D_{x_2}U_2^n=-D_{x_1}U_1^n\,,\mbox{ for all }(t,x)\mbox{ and }n\in\N\,. $$ Since for all $t>0$, the kernels are of Calderon-Zigmund kind, via \rf{BN}, we can deduce uniformly in $t>0$ and in $n\in\N$ \be\label{IP-I}\ba{ll}\dm D_{x_1}\!U_1^n(t)\dm_p\!+ \!\dm D_{x_1}\!U_2^n(t)\dm_p\!+\!\dm D_{x_2}\!U_2^n(t)\dm_p\hskip-0.3cm&\leq \!c\dm \Gamma[\nabla u_0^n](t)\dm_p\!\leq\! c\dm \n u^n_0\dm_p\VSE\leq\! c\dm\n u_0\dm_p+o(1)\,,\ea\ee where $\Gamma[\n u_0^n](t,x):= \intl{\R^2_+}\Gamma(t,x-y)\n u_0(y)dy\,.$
   We estimate $D_{x_2}U^n_1$. Since $\frak G_{12}=0$, we restrict ourselves to consider only $\frak G_{11}[{u_0^n}_1]$. Since $u_0^n$ has compact support \be\label{WQP}\mbox{for all } q\in[1,\infty],\quad\Gamma[u_0^n],\Gamma^*[u_0^n]\in W^{1,q}(\R^2_+),\mbox{ for all }t>0.\ee Hence, integrating by parts with respect to $y_1$, we easily get
$$\ba{l}\displ U^n_1(t,x)=   (\Gamma -\Gamma^*)[{u^n_0}_1](t,x)
 +4\!\!\intll0{x_2}\!\!\intl{\R}\!\!D_{x_1}E(x\!-\!z)\Gamma^*[{D_{y_1}u^n_0}_1](t,z)dz=I_1^n+I_2^n\,.\ea$$
Integrating by parts with respect to $y_2$, we get $$D_{x_2}I^n_1(t,x)=(\Gamma+\Gamma^*)[D_{y_2}{u_0^n}_1]\,.$$ Hence, it follows
\be\label{IP-II}\dm D_{x_2}I_1^n(t)\dm_p\leq c\dm \n u_0\dm_p+O(n^{-1+\frac2p})\,\mbox{ for all }t>0\,\mbox{ and }n\in\N\,.\ee
In the case of $I_2^n$, by virtue of \rf{WQP},  we get
$$\ba{l}\displ D_{x_2}I^n_2=-4\,{\rm P.V.}\hskip-0.1cm\intll0{x_2}\intl{\R}D^2_{z_2z_1}E(x-z)\Gamma^*[D_{y_1}u_0^n](t,z)dz  \VS\hskip4.5cm-4\,{\rm P.V.}\hskip-0.1cm\intl\R D_{z_1}E(x_1-z_1,0)\Gamma^*[D_{y_1}u_0^n](t,z_1,x_2)dz_1\,.  \ea$$ Applying the Calderon-Zigmund theorem and the properties of heat kernel, we deduce the estimate \be\label{IP-III}\dm D_{x_2}I_2^n(t)\dm_p\leq c\dm \n u_0^n\dm_p\leq c\dm \n u_0\dm_p+o(1)\,\mbox{ for all }t>0\mbox{ and }n\in\N.\ee Collecting estimates \rf{IP-I} and \rf{IP-II}-\rf{IP-III}, we deduce that
$$\dm \n U^n(t)\dm_p\leq c\dm \n u_0\dm_p+o(1), \mbox{ for all }t>0\mbox{ and }n\in\N\,.$$
This last estimate, the pointwise convergence ensured by \rf{PGE} and the Fatou theorem prove \rf{GPST-I}. If $p\in(1,2)$, we get $L^\infty(\R^2_+)\cap J_0^p(\R^2_+)\subset L^\infty(\R^2_+)\cap L^q(\R^2_+)$, provided that $q=\frac{2p}{2-p}$\,. Hence all the the above computations hold without  the special approximation $\{u_0^n\}$ of Lemma\,\ref{DEN}.
 \ep We conclude the section recalling the following well known and special result (see e.g. \cite{MS-II}):\begin{tho}\label{LPSP}{\sl Let $u_0=0$ in \rf{ST}. For all $G\in L^r(0,T;L^r(\R^2_+))$, $r\in(1,\infty)$, there exists a unique solution  to problem \rf{ST} such that $U\in C(0,T;J^{1,r}(\R^2_+))$ with \be\label{LPSP-I} \ba{c}\displ t^{-\frac{r-1}r}\dm U(t)\dm_r+t^{-\frac{r-2}{2r}}\dm \n U(t)\dm_r\!\leq c \Big[\intll0t\!\dm G(\tau)\dm_r^rd\tau\Big]^\frac1r\!,\mbox{\,a.e.\,in\,}t\!\in\!(0,T),\\\displ\intll0T\Big[\dm U_t(t)\dm_r^r+\dm D^2U(t)\dm_r^r+\dm \n\pi_U(t)\dm_r^r\Big]dt\leq c\intll0T\dm G(t)\dm_r^rdt\,,\ea\ee where the constant $c$ is independent of   $G$. If $u_0\in \mathscr C_0(\R^2_+)$, then \rf{LPSP-I}$_2$ becomes
\be\label{LEPSP}\displ\intll0T\!\Big[\dm U_t(t)\dm_r^r+\dm D^2U(t)\dm_r^r+\dm \n\pi_U(t)\dm_r^r\Big]dt\leq c\Big[\dm u_0\dm_{{2-\frac2r,r}}^r+\!\intll0T\!\dm G(t)\dm_r^rdt\Big],\ee with $c$ indepedndent of $u_0$ and $G$.
 }\end{tho}
\section{\label{SLNS}A linearezed Navier-Stokes problem}  
For our aims   the following initial boundary value problem is crucial:
\be\label{LNS}\ba{l}v_t\!-\Delta v+\!\n \pi_v\!=-\mbox{${\overset{h}{\underset{\ell=0}\sum}}$}w^\ell\!\cdot\!\!\n v-v\!\cdot\!\!\n \mbox{${\overset{h}{\underset{\ell=0}\sum}}$}w^\ell\!+\!F\,,\VS\n\!\cdot\! v=0\mbox{\quad\,in\,\;}(0,T)\!\!\times\!\R^2_+,\VS v=0\mbox{\quad on }(0,T)\times\{x_2=0\},\quad v=0\mbox{\quad on }\{0\}\times\R^2_+.\ea\ee
\begin{ass}{\sl For problem \rf{LNS} we assume:\begin{itemize}\item[i.] $w^0\in L^\infty((0,\infty)\times \R^2_+)\cap L^\infty(0,\infty;J^p_0 (\R^2_+))\,,$ and $t^\frac12| \nabla w^0(t)|_\infty\leq c\dm w^0(0)\dm_\infty$\,.
\item[ii.] For all $\ell=1,\dots,h$, and for all $T>0$,  \newline$w^\ell\in C(0,T;J^{1,\frac p{2^{\ell-1\hskip-0.3cm\null}}}\hskip0.3cm(\R^2_+))\cap L^{ \frac p{2^{\ell-1\hskip-0.3cm\null}}}\hskip0.25cm(0,T;W^{2,{\frac p{2^{\ell-1\hskip-0.3cm\null}}}}\hskip0.25cm(\R^2_+))$\,, where we assume $p>2^{\ell+1}$\,.\item[iii.]$F\in L^{\frac p{2^h\hskip-0.1cm}}\,(0,T;L^{\frac p{2^h\hskip-0.1cm\null}}\,(\R^2_+))$, where we assume $\frac p{2^h\hskip-0.1cm\null}>2$\,.\end{itemize}}\end{ass}We start with the following
\begin{lemma}\label{L-I}{\sl Assume\, {\rm i.-ii.} for $w^\ell$, $\ell=0,\dots,h$. Moreover assume that $v_t,D^2v\in L^{\frac p{2^h\hskip-0.1cm\null}}\,(0,T; L^{\frac p{2^h\hskip-0.1cm\null}}\,(\R^2_+))$ with $v=0$ a.e. on $(0,T)\times\{x_2=0\}$ and $v=0$ a.e. on $\{0\}\times \R^2_+$. Then, a.e. in $t\in(0,T)$, it holds that
\be\label{L-I-I}\ba{l}\dm w^\ell(t)\cdot\n v(t)\dm_{\frac p{2^h\!\!\!\null}}\leq c\dm w^\ell(t)\dm_\infty\dm D^2v(t)\dm_{\frac p{2^h\!\!\!\null}}^\frac12\,\Big[\intll0t\dm v_\tau(\tau)\dm_{\frac p{2^h\!\!\!\null}}d\tau\,\Big]^\frac12\;,\\\dm v(t)\cdot\n w^\ell(t)\dm_{\frac p{2^h\!\!\!\null}}\leq c \dm\n w^\ell(t)\dm_\infty\intll0t\dm v_\tau(\tau)\dm_{\frac p{2^h\!\!\!\null}}\,d\tau \,,  \ea\ee for $\ell=0,\dots,h$, and with $c$ independent of $w^\ell,\,v$ and $t$.}   \end{lemma}\bp Sine $p/2^{\ell-1}>2$, empolying the Gagliardo-Nirenberg inequality, for all $\ell=1,\dots,h$, we get ($q_\ell:=\frac p{2^{\ell-1}\!\!\!\!\!\!\!\!\null}$\hskip0.35cm)
\be\label{GNW}\ba{l}\dm w^\ell\dm_\infty\leq c\dm D^2w^\ell\dm_{q_\ell}^a\;\;\dm w^\ell\dm_{q_\ell}^{1-a}<\infty\mbox{ a.e. in }t>0\,,\,a=\frac1{q_\ell\!\!\!\null}\,,\VS\dm \n w^\ell\dm_\infty\leq c\dm D^2w^\ell\dm_{q_\ell}^b\dm w^\ell\dm_{q_\ell}^{1-b}<\infty\mbox{ a.e. in }t>0\,,\mbox{$b=\frac12+\frac1{q_\ell\!\!\!\null}  $}\,.\ea\ee We   recall  that in our hypotheses on $v$ the following estimates hold
\be\label{L-I-II}\dm v(t)\dm_{\frac p{2^h\!\!\!\null}}\leq \intll0t\dm v_\tau(\tau)\dm_{\frac p{2^h\!\!\!\null}}\;d\tau\,\mbox{ for all  }t\in(0,T)\,,\ee and \be\label{L-I-III}\! \dm \n v(t)\dm_{\frac p{2^h\!\!\!\null}}\!\leq c\dm D^2 v(t)\dm_{\frac p{2^h\!\!\!\null}}^\frac12\dm v(t)\dm_{\frac p{2^h\!\!\!\null}}^\frac12\,,\mbox{ a.e. in }t\in(0,T)\,.  \ee
Applying Holder's inequality, after employing estimates \rf{L-I-II}-\rf{L-I-III}, estimates \rf{L-I-I} easily follow.\ep
\begin{lemma}\label{L-O}{\sl Assume\, {\rm i.-ii.} for $w^\ell$, $\ell=0,\dots,h$. Moreover assume that $v\in C(0,T;J^{1,2}(\R^2_+))$. Then, a.e. in $t\in(0,T)$ the following estimates hold with $q:={\frac p{2^h\!\!\!\!\null}}\hskip0.2cm$
\be\label{L-O-I}\ba{l}\dm w^\ell(t)\cdot\n J_\vep[ v(t)]\dm_q\leq c(\vep)\dm w^\ell(t)\dm_\infty\dm \n v(t)\dm_2 \;,\VS\dm J_\vep[v(t)]\cdot\n w^\ell(t)\dm_q\leq c(\vep) \dm  v(t)\dm_2\dm\n w^\ell(t)\dm_\infty \,,\,  \ea\ee for all $\ell=0,\dots,h$, here $J_\vep[\cdot] $ is a spatial mollifier   and constant $c(\vep)$  is independent of $w^\ell,\,v$ and $t$.}   \end{lemma}\bp Estimates \rf{L-I-I} hold. Since $q>2$, by virtue of the properties of the mollifier the inequalities \rf{L-O-I} hold immediately.  \ep
\begin{lemma}\label{LD-I}{\sl Assume {\rm i.-ii.} for $w^\ell,\,\ell=0,\dots,h$. Moreover, assume that $v\in L^2(0,T;J^{1,2}(\R^2_+))$ and $V\in L^2(0,T;L^2(\R^2_+))$. Then, a.e. in $t>0$, for $p>2^\ell$, the following estimates hold
\be\label{LD-II}\ba{l}|(w^0(t)\cdot\n J_\vep[v(t)],V(t))|\leq \dm w^0(t)\dm_\infty\dm \n v(t)\dm_2\dm V(t)\dm_2\,,\VS |(w^\ell(t)\cdot\n J_\vep[v(t)],V(t))|\leq c(\vep) \dm w^\ell(t)\dm_{\frac p{2^{\ell-1}\hskip-0.3cm\null}}\hskip0.3cm\dm \n v(t)\dm_2\dm V(t)\dm_2,\VS |(J_\vep[v(t)]\cdot\n w^0(t),V(t))|\leq c(\vep)\dm \n w^0(t)\dm_{p }\dm v(t)\dm_2\dm V(t)\dm_2\,,\VS |(J_\vep[v(t)]\cdot\n  w^\ell(t),V(t))|\leq c(\vep)\dm\n w^\ell(t)\dm_{\frac p{2^{\ell-1\hskip-0.3cm\null}}}\hskip0.3cm\dm  v(t)\dm_2\dm V(t)\dm_2\,,\ea\ee
where $J_\vep[\cdot]$ is a spatial mollifier.
}\end{lemma}\bp Applying Holder's inequality we get$$\ba{l}|(w^0\cdot\n J_\vep[v],V)|\leq \dm w^0\dm_\infty\dm \n v\dm_2\dm V\dm_2\,,\VS |(w^\ell\cdot\n J_\vep[v],V)|\leq  \dm w^\ell\dm_{\frac p{2^{\ell-1\hskip-0.3cm\null}}}\hskip0.4cm\dm \n J_\vep[ v]\dm_{\frac{2p}{p-2^\ell\hskip-0.1cm}}\hskip0.3cm\dm V\dm_2\,,\VS |(J_\vep[v]\cdot\n w^0,V)|\leq  \dm \n w^0\dm_{p }\dm J_\vep[ v]\dm_{\frac{2p}{p-2}}\dm V\dm_2\,,\VS |(J_\vep[v]\cdot\n  w^\ell,V)|\leq  \dm\n w^\ell\dm_{\frac p{2^{\ell-1}\hskip-0.3cm\null}}\hskip0.3cm\dm J_\vep[ v]\dm_{\frac{2p}{p-2^\ell\hskip-0.1cm}}\hskip0.1cm\dm V\dm_2\,.\ea$$Hence, by virtue of properties of the mollifier, and our hypotheses on $v$ and $V$, we deduce \rf{LD-II} a.e. in $t\in(0,T)$.\ep
 \begin{tho}\label{TLNS}{\sl Under Assumption\,1 for $w^\ell$ and $F$, there exists a unique solution to problem \rf{LNS} such that, for all $T>0$, $v\in C(0,T;J^{1,q}(\R^2_+))$ with $q:={\frac p{2^h\!\!\!\!\null}}$\hskip0.2cm \be\label{TLNS-I}\ba{c}\displ t^{-\frac{q-1}q}\dm v(t)\dm_{q}\!+t^{-\frac{q-2}{2q}}\dm \n v(t)\dm_q\! \leq \!c\,{\rm exp}[\intll0t\!g(t)dt]\Big[\intll0t\!\dm F(t)\dm_q^qdt\Big]^\frac1q,\\\hskip9cm\mbox{\,for all }t\!\in\![0,T)\,,\\\displ\intll0T\!\Big[\dm v_t(t)\dm_q^q\!+\dm D^2v(t)\dm_q^q\!+\dm \n\pi_v(t)\dm_q^q\Big]\!\leq c\,{\rm exp}[\intll0T\!\!cg(t)dt]\!\intll0T\!\!\dm F(t)\dm_q^qdt \,,\ea\ee where the constant $c$ is independent of $w^\ell$ and $F$, and we set $g(t):=t^{q-1}\,\mbox{${\overset{h}{\underset{\ell=0}\sum}}$}\big[\dm w^\ell(t)\dm_{ \infty}^{2q}+\dm\n w^\ell(t)\dm_\infty^q\big]$.}\end{tho} 
\bp We introduce the following approximation  problem:
\be\label{ALNS}\ba{c}v_t-\Delta v+\!\n \pi_v=-\mbox{${\overset{h}{\underset{\ell=0}\sum}}$}w^\ell\!\cdot\!\n J_\vep[v]-J_\vep[v]\cdot\!\n \mbox{${\overset{h}{\underset{\ell=0}\sum}}$}w^\ell+ F_\vep\,,\VS\n\!\cdot\! v=0\,,\mbox{ in }(0,T)\!\times\!\R^2_+,\VS v=0\mbox{ on }(0,T)\times\{x_2=0\},\quad v=0\mbox{ on }\{0\}\times\R^2_+,\ea\ee
where $J_\vep[\cdot]$ is a spatial   mollifier and $F_\vep:={\rm exp}[-\vep|x|^2]F$. Of course, a solution to problem \rf{ALNS} is a pair $(v_\vep,\pi_{v_\vep})$, for the sake of the simplicity and of the brevity we omit the index $\vep$. For all $\vep>0$ and $T>0$, the data $F_\vep$ 
belongs to $L^2(0,T;L^2(\R^2_+))$. Thanks to this firstly  we are able to develop  a $L^2$-theory for problem \rf{ALNS} depending on $\vep$. Then on the base of this result, we approach   the $L^q$-theory of solution to problem \rf{LNS}  based on the $L^q$-theory of the Stokes problem. Employing the Galerkin method, in the way suggested by Heywood in \cite{H},  we can easily establish the existence of  the Galerkin approximation sequence which satisfies the set of relations
$$\mbox{$\frac12$}\dm v\dm_2^2+\intll0t\dm\n v\dm_2^2d\tau =-\intll0t(\mbox{${\overset{h}{\underset{\ell=0}\sum}}$}w^\ell\!\cdot\!\n J_\vep[v]+J_\vep[v]\cdot\!\n \mbox{${\overset{h}{\underset{\ell=0}\sum}}$}w^\ell,v)d\tau+\intll0t(F_\vep,v)d\tau\,,$$
$$\mbox{$\frac12$}\dm\n v\dm_2^2+\intll0t\dm P\Delta v\dm_2^2d\tau =\intll0t(\mbox{${\overset{h}{\underset{\ell=0}\sum}}$}w^\ell\!\cdot\!\n J_\vep[ v]+J_\vep[v]\cdot\!\n \mbox{${\overset{h}{\underset{\ell=0}\sum}}$}w^\ell,P\Delta v)d\tau-\intll0t(F_\vep,P\Delta v)d\tau\,,$$
$$\mbox{$\frac12$}\dm\n v\dm_2^2+\intll0t\dm  v_t\dm_2^2d\tau =-\intll0t(\mbox{${\overset{h}{\underset{\ell=0}\sum}}$}w^\ell\!\cdot\!\n J_\vep[ v]+J_\vep[v]\cdot\!\n \mbox{${\overset{h}{\underset{\ell=0}\sum}}$}w^\ell,  v_t)d\tau+\intll0t(F_\vep,v_t)d\tau\,.$$ Applying H\"older's inequality, employing estimates \rf{LD-II} with $V$ substituted by $v$, and by   Lemma\,\ref{GW} from the first relation of the set we easily get the energy inequality
\be\label{EI}\ba{l}\displ\dm v(t)\dm_2^2+\!\intll0t\dm\n v(\tau)\dm_2^2d\tau\leq c\,  {\rm exp}\Big[t+c(\vep)\!\intll0t\!\!\big[\dm w^0\dm^2_\infty\!+\mbox{${\overset h{\underset {\ell=1}\sum}}$}\dm w^\ell\dm_{\frac p{2^{\ell-1}\!\!\!\!\!\!\!\!\!\null}}^2 \hskip0.4cm\big]d\tau \VS\hskip2.4cm+c(\vep)\! \intll0t\!\!\big[\dm\n w^0\dm_p\!+\!\mbox{${\overset h{\underset {\ell=1}\sum}}$}\dm\n w^\ell\dm_{\frac p{2^{\ell-1}\!\!\!\!\!\!\!\!\!\null}}\,\;\;\big]d\tau\! \!\intll0t\!\dm F_\vep\dm_2^2\;d\tau=:\!\frak B(t).\ea\ee Subsequently, from the remaining two relations of the above set,   applying H\"older's inequality, employing estimates \rf{LD-II} with $V$ substituted by $P\Delta v,\,v_t$ along the two cases, and via \rf{EI}   we easily get ,  $$\ba{l}\displ\dm \n v(t)\dm_2^2+\intll0t\big(\dm v_\tau(\tau)\dm_2^2+\dm P\Delta v(\tau)\dm_2^2\big)d\tau\leq  c(\vep)\Big[t\frak B(t) \sup_{(0,T)}\big(\dm w^0\dm^4_\infty+\dm\n w^0\dm_p^2 \VS\hskip5.8cm+\mbox{${\overset h{\underset {\ell=1}\sum}}$}\dm w^\ell\dm_{\frac p{2^{\ell-1}\!\!\!\!\!\!\!\!\!\null}}^4\;\;+ \mbox{${\overset h{\underset {\ell=1}\sum}}$}\dm \n w^\ell\dm_{\frac p{2^{\ell-1}\!\!\!\!\!\!\!\!\!\null}}^2\;\;\;\big)+\intll0t\dm F_\vep\dm_2^2d\tau\Big]. \ea$$ Now standard arguments related to the Galerkin method allow  us to deduce the existence of a pair $(v_\vep,\pi_{v_\vep})$ solution to problem \rf{ALNS} such that 
\be\label{LDE}\ba{l}v_\vep\! \in\! C([0,T,J^{1,2}(\R^2_+))\cap L^2(0,T;W^{2,2}(\R^2_+)),\VS{v_\vep}_t\,,\!\n\pi_{v_\vep}\!\!\in\! L^2(0,T;L^2(\R^2_+))\,.\ea\ee
Now our task is to prove that the family of solutions $\{v_\vep\}$ admits a limit for $\vep\to0$ enjoying of the property \rf{TLNS-I}. By virtue of \rf{LDE}, we can apply Lemma\,\ref{L-O} for each term of the right hand side of problem \rf{ALNS}. Hence, for all $\vep>0$, we can claim that 
$$-\mbox{${\overset{h}{\underset{\ell=0}\sum}}$}w^\ell\!\cdot\!\n J_\vep[v]-J_\vep[v]\cdot\!\n \mbox{${\overset{h}{\underset{\ell=0}\sum}}$}w^\ell+ F_\vep\in L^q(0,T;L^q(\R^2_+))\,.$$
By virtue of Theorem\,\ref{LPSP} related to the Stokes problem, for all $\vep>0$ and $T>0$, we can claim that
\be \label{LPASP-I} \ba{c}s^{-\frac{q-1}q}\displ\dm v(s)\dm_q+s^{-\frac{q-2}{2q}}\dm \n v(s)\dm_q\leq c \Big[\intll0s\dm G(\tau)\dm_q^qd\tau\Big]^\frac1q,\mbox{ for all }s\in(0,T)\,,\\\displ\intll0s\!\!\Big[\dm v_t(t)\dm_q^q+\dm D^2v(t)\dm_q^q+\dm \n\pi_v(t)\dm_q^q\Big]dt\!\leq c\!\!\intll0s\dm G(t)\dm_q^qdt\,,\mbox{ for all }s\!\in\!(0,T)\ea\ee
where $G:=-\mbox{${\overset{h}{\underset{\ell=0}\sum}}$}w^\ell\!\cdot\!\n J_\vep[v]-J_\vep[v]\cdot\!\n \mbox{${\overset{h}{\underset{\ell=0}\sum}}$}w^\ell+ F_\vep\in L^q(0,T;L^q(\R^2_+))$\,, and $c$ is independent of $\vep$ and $T$. Now we look for estimates in $L^q(0,T;L^q(\R^2_+))$ for function $G$ which are independent of $\vep$. By virtue of Lemma\,\ref{L-I} and estimates \rf{GNW} for $w^\ell$, applying the Cauchy and the H\"older inequality, we get  
\be\label{IDE-I}  \dm w^\ell(t)\!\cdot\!\n\! J_\vep[v(t)]\dm_q^q \leq c t^{q-1}\;\dm w^\ell(t)\dm_\infty^{2q}\intll0t\dm v_\tau(\tau)\dm_q^qd\tau +\eta\dm D^2 v(t)\dm_q^q \,,\ee   for all $\ell=0,\dots,h$. Analogously, we get
\be\label{IDE-II} \dm J_\vep[v(t)]\cdot\n w^\ell(t)\dm_q^q \leq ct^{{q-1}}\,\dm \n w^\ell(t)\dm_\infty^q\intll0t\dm v_\tau(\tau)\dm_q^qd\tau \,,\ee for all $\ell=0,\dots,h$. Finally, we easily deduce that \be\label{FP}\dm F_\vep(t)\dm_q\leq \dm F(t)\dm_q\,,\mbox{ for all }t>0\mbox{ and }\vep>0\,.\ee Collecting estimates \rf{IDE-I}-\rf{FP}, recalling the definition of $G$, from \rf{LPASP-I}$_2$, for a suitable $\eta>0$ in estimate \rf{IDE-I}, for all $s>0$, we obtain
\be\label{IDE-III}\hskip-0.2cm\intll0s\!\!\big[\dm v_t(t)\dm_q^q\!+\dm D^2v(t)\dm_q^q\!+\dm \n\pi_v(t)\dm_q^q\big]dt \!\leq\! c\!\!  \intll0s\!\big[g(t) \!\! \intll0t\!\!\dm v_\tau(\tau)\dm_q^qd\tau\!+\dm F(t)\dm_q^q\big]dt, \ee where we have set $$g(t):=t^{q-1}\,\mbox{${\overset{h}{\underset{\ell=0}\sum}}$}\big[\dm w^\ell(t)\dm_{ \infty}^{2q}+\dm\n w^\ell(t)\dm_\infty^q\big]\,.$$ From \rf{IDE-III} an application of   Lemma\,\ref{GW} ensures that
\be\label{IDE-V}\intll0s\dm v_t(t)\dm_q^qdt\leq c\,{\rm exp}[\intll0sg(t)dt]\intll0s\dm F(t)\dm_q^qdt\,.\ee Enclosing the last estimate in the right hand side of \rf{IDE-III} a trivial computation furnishes
\be\label{IDE-IV}\ba{l} \displ\intll0s\!\!\big[\dm v_t(t)\dm_q^q\!+\dm D^2v(t)\dm_q^q\!+\dm \n\pi_v(t)\dm_q^q\big]dt\leq  c\,{\rm exp}[\intll0s\!\!g(t)dt]\!\intll0s\!\!\dm F(t)\dm_q^qdt\,.   \ea\ee
Estimates \rf{IDE-V}-\rf{IDE-IV} are independent of $\vep$. Hence, taking into account the definition of $G$, collecting estimates \rf{IDE-I}-\rf{FP} and estimate \rf{IDE-IV},    we have proved that $$\intll0s\dm G(\tau)\dm_q^qd\tau\leq c\exp[\intll0s g(t)dt]\intll0s\dm F(t)\dm_q^qdt\,, \;s\in(0,T)\,.$$ Thus,   proving that $\intll0s g(t)dt$ is finite for all $s>0$, then we have concluded the proof of the theorem. By virtue of the assumption on $w^0$ and weight $t^{q-1}$, the  integrability question can be restricted to the cases of $\ell=1,\dots,h$. By virtue of estimates \rf{GNW}$_1$ and assumptions on $w^\ell$, we have
$$\dm w^\ell\dm_\infty^{2q}\leq  c\dm D^2w^\ell\dm^{A_\ell}_{\frac p{2^{\ell-1\!\!\!\!\!\!\!\!\!\null}}}\hskip0.35cm\dm w^\ell\dm_{\frac p{2^{\ell-1\!\!\!\!\!\!\!\!\!\null}}}^{B_\ell}\hskip0.2cm\leq c\sup_{(0,T)}\dm w^\ell\dm_{\frac p{2^{\ell-1\!\!\!\!\!\!\!\!\!\null}}}^{B_\ell}\hskip0.25cm \dm D^2w^\ell\dm^{A_\ell}_{\frac p{2^{\ell-1\!\!\!\!\!\!\!\!\null}}}\hskip0.3cm,$$ where $A_\ell:=2qa_\ell={2^{\ell-h}}$ and $B_\ell:=2q(1-a_\ell)={2^{1-h}(p-2^{\ell-1})}$.  Since $p>2^{\ell+1}$, we obtain $$\intll0s\dm w^\ell\dm_\infty^{2q}dt\leq c\sup_{(0,T)}\dm w^\ell\dm_{\frac p{2^{\ell-1\!\!\!\!\!\!\!\!\!\null}}}^{B_\ell}\;\;\intll0s \dm D^2w^\ell\dm^{A_\ell}_{\frac p{2^{\ell-1\!\!\!\!\!\!\!\!\!\null}}}dt<\infty \mbox{ for all }s>0 \mbox{ and }\ell=1,\dots,h\,.$$ Analogously, by virtue of estimates \rf{GNW}$_2$, we get
$$\dm \n w^\ell\dm_\infty^q\leq c\dm w^\ell\dm_{\frac p{2^{\ell-1\!\!\!\!\!\!\!\!\!\null}}}^{C_\ell}\hskip0.3cm\dm D^2w^\ell\dm^{D_\ell}_{\frac p{2^{\ell-1\!\!\!\!\!\!\!\!\!\null}}}\hskip0.2cm\leq c\sup_{(0,T)}\dm w^\ell\dm_{\frac p{2^{\ell-1\!\!\!\!\!\!\!\!\!\null}}}^{C_\ell}\hskip0.25cm\dm D^2w^\ell\dm^{D_\ell}_{\frac p{2^{\ell-1\!\!\!\!\!\!\!\!\!\null}}}$$  where $C_\ell:=q(1-b_\ell)= \frac p{2^h\!\!\!\null}\;(\frac12-\frac{2^{\ell-1\!\!\!\!\!\!\!\!\!\null}}p\;\;\;)$ and $D_\ell:=qb_\ell=\frac p{2^h\!\!\!\null}\;(\frac 12+\frac{2^{\ell-1\!\!\!\!\!\!\!\!\!\null}}p\;\;\;)$. Since $p>2^{\ell+1}$ we obtain 
$$\intll0s\!\!\dm\n w^\ell\dm_\infty^{q}dt\leq c\sup_{(0,T)}\dm w^\ell\dm_{\frac p{2^{\ell-1\!\!\!\!\!\!\!\!\!\null}}}^{C_\ell}\;\,\,\intll0s \!\!\dm D^2w^\ell\dm^{D_\ell}_{\frac p{2^{\ell-1\!\!\!\!\!\!\!\!\!\null}}}dt<\infty \mbox{ for all }s>0 \mbox{ and }\ell=1,\dots,h\,.$$ The theorem is completly proved.
\ep  
\section{\label{PU}The linearezed Navier-Stokes IBVP in 
 \mbox{\;\footnotesize$\bf\cap_{_{_{\hskip-0.5cm1<q<2}}} J^q$}}
 In order to discuss the uniqueness we have to consider the   following linearezed Navier-Stokes problem:
 \be\label{LNSQ}\ba{l}\varphi_t-\Delta \varphi+\!\n \pi_\varphi\!=V\!\cdot\!\n \varphi+\varphi\!\cdot\!(\n\mbox{$\overset \nu{\underset{h=1}\sum}$} W^h)^T,\;\nabla\!\cdot\!\varphi=0\mbox{ in }(0,T)\!\times\!\R^2_+,\VS \varphi=0\mbox{ on }(0,T)\times\{x_2=0\},\quad \varphi=\varphi_\circ\mbox{ on }\{0\}\times\R^2_+,\ea\ee
here the symbol $(\n b)^T$ means the transpost of tensor $\n b$ and $a\cdot(\n b)^T= (\frac{\partial b_i}{\partial x_k}a_i)e_k$. We assume that for all $T>0$\be\label{LNSQO} V\!\in\! C((0,T)\times\R^2_+),\mbox{ and }\n W^h\!\in\! C(0,T;L^{r^h}(\R^2_+))\mbox{\,for\,some\,}r^h>2\,.\ee The investigation on problem \rf{LNSQ} appears very similar to the one of the previous section. Actually there are some different technical aspects that make the difference. Theorem\,\ref{TLNS} is related to the existence, whereas Theorem\,\ref{TLNSQ}   given below is related to   the uniqueness by duality. Although the proofs are  approached by a similar way (that is, initially by means of the $L^2$-theory), we look on the two theorems from a different point of view. As opposed  to the previous section, here we discuss the initial boundary value problem with an initial data $\varphi_\circ\ne 0$ and body force $F=0$,   we study an $L^q$ theory for $q\in(1,2)$ with the special property \rf{LNSQ-I} (see below Theorem\,\ref{TLNSQ}). These are not given in section\,\ref{SLNS}. As well the following lemmas  are thought  by slightly different way. Hence, in order to make more readable the results of section\,\ref{SLNS} and of this section, we have thought that it is better to furnish the results in two separated theorems, rather than to state all the results in a unique large theorem.\par     We start with the following 
\begin{lemma}\label{LN-I}{\sl Let $\varphi\in J^{1,2}(\R^2_+)$ and let $\phi\in L^2(\R^2_+)$. Then it holds that
\be\label{L-II}\ba{c}|(V\cdot J_\vep[\chi_\vep\n\varphi],\phi)| \leq \dm V\dm_\infty\dm\n\varphi\dm_2\dm\phi\dm_2 \VS |(J_\vep[\varphi]\cdot\!(\n W^h)^T,\phi)| \leq c \dm \n\varphi\dm_2^{a_h}\dm \varphi\dm_2^{1-a_h}\dm \n W^h\dm_{r^h}\dm \phi\dm_2 \,,\ea   \ee with $a_h=\frac2{r^h\!\!\null}$\,\, and $c$ independent of $\varphi,\,\phi\,.$ Here $J_\vep[\cdot]$ is  a spatial mollifier, and $\chi_\vep(x):=\chi(\vep x)$ is the smooth  cutoff function with support $(\frac1\vep,\frac2\vep)$ ($\chi$ is defined in section\,2).  
 }\end{lemma}\bp Applying H\"older's inequality, employing the properties of the mollifier, and Gagliardo-Nirenberg inequality we get \be\label{LNSQ-I}\ba{ll}\null\hskip0.5cm|(V\cdot J_\vep[\chi_\vep\n\varphi],\phi) |  \hskip-0.2cm&\leq \dm V\dm_\infty\dm\n\varphi\dm_2\dm\phi\dm_2 \VS   |(J_\vep[\varphi]\cdot\!(\n W)^T,\phi)|&\leq\dm \varphi\dm_{\frac{2r}{r-2}}\dm \n W\dm_{r}\dm \phi\dm_2\VS&\leq c \dm \n\varphi\dm_2^a\dm \varphi\dm_2^{1-a}\dm \n W\dm_{r}\dm \phi\dm_2 \,,\ea\ee which prove \rf{L-II}\ep\begin{lemma}\label{UESP}{\sl    Let $q\leq r_h$ and $\varphi\in C(0,T;J^q(\R^2_+))\cap L^q(0,T; J^{1,q}(\R^2_+)\cap W^{2,q}(\R^2_+))$. We assume that $\varphi_t\in L^q(0,T;L^q(\R^2_+))$. Then almost everywhere in $t\in(0,T)$ it holds that 
\be\label{UESP-I}\ba{c}\displ\dm V\cdot J_\vep[\chi_\vep\n\varphi]\dm_q \leq c\dm V\dm_\infty\dm D^2\varphi\dm_q^\frac12\Big[\intll0t\dm\varphi_\tau\dm_qd\tau+\dm \varphi_\circ\dm_q\Big]^\frac12 \VS \dm J_\vep[\varphi]\cdot\!(\n W^h)^T\dm_q \leq c  \dm \n W^h\dm_{r^h}\dm D^2\varphi\dm_q^\frac1{r_h\!\!\!\null}\,\Big[\intll0t\dm \varphi_\tau\dm_qd\tau+\dm \varphi_\circ\dm_q\Big]^{\frac1{r'_h\!\!\!\null}} \,\,,\ea   \ee  where $J_\vep[\cdot]$ is  a spatial mollifier, and $\chi_\vep(x):=\chi(\vep x)$ is the smooth  cutoff function with support $(\frac1\vep,\frac2\vep)$ ($\chi$ is defined in section\,2). The constant $c$ is independent of $\varphi$ and $\vep$. }\end{lemma}  \bp We recall the following: 
$$\dm \varphi(t)\dm_q\leq \dm \varphi_\circ\dm_q+\intll0t\dm \varphi_\tau\dm_qd\tau\,,\mbox{ for all }t\in(0,T)\,.$$Applying H\"older's inequality and Gagliardo-Nirenberg's inequality, we get
$$ \dm V\cdot J_\vep[\chi_\vep\n\varphi]\dm_q \leq c\dm V\dm_\infty\dm D^2\varphi\dm_q^\frac12\dm \varphi\dm_q^\frac12\,,\mbox{ a.e. }t\in(0,T).$$
So that easily estimate \rf{UESP-I}$_1$ follows. Analogously applying H\"older's inequality and Gagliardo-Nirenberg's inequality we get
$$\dm J_\vep[\varphi]\cdot(\n W)^T\dm_q \leq \dm \n W\dm_r\dm \varphi\dm_{\frac{rq}{r-q}}\leq c\dm \n W\dm_r\dm D^2\varphi\dm_q^\frac1r\dm \varphi\dm_q^{\frac1{r'\!\!\null}}\,,\mbox{ a.e. in }t\in(0,T)\,.$$ Again we claim that easily estimate \rf{UESP-I}$_2$ follows.
\ep
\begin{tho}\label{TLNSQ}{\sl  For all $\varphi_\circ\in \mathscr C_0(\R^2_+)$ there exists a unique solution to problem \rf{LNSQ} such that for all $T>0$ \be\label{LNSQ-I}\ba{c}\varphi\in {\underset{1<q<2}\cap}\Big[C(0,T;J^q(\R^2_+))\cap L^q(0,T;J^{1,q}(\R^2_+))\Big],\VS D^2\varphi,\,\varphi_t,\n\pi_\varphi\in {\underset{1<q<2}\cap}L^q(0,T;L^q(\R^2_+))\,.\ea\ee}\end{tho}\bp   We consider the approximation problem
\be\label{LNSQ-II}\ba{l}\varphi_t-\Delta \varphi+\!\n \pi_\varphi=V\!\cdot\! J_\vep[\chi_\vep\n\varphi]+J_\vep[\varphi]\cdot\!(\n{\underset{h=1}{\overset\nu\sum}} W^h)^T,\\\nabla\cdot\varphi=0\mbox{ in }(0,T)\!\times\!\R^2_+,\VS \varphi=0\mbox{ on }(0,T)\times\{x_2=0\},\quad \varphi=\varphi_\circ\mbox{ on }\{0\}\times\R^2_+,\ea\ee where $J_\vep[\cdot]$ is  a spatial mollifier, and $\chi_\vep(x):=\chi(\vep x)$ is the smooth  cutoff function with support $(\frac1\vep,\frac2\vep)$ ($\chi$ is defined in section\,2).  In order to obtain the solution to problem \rf{LNSQ-II}, we can apply as usual the Galerkin method.  Employing estimates \rf{L-II}, where in place of $\phi$ we set $\varphi$ to obtain the energy inequality, then $\phi=P\Delta \varphi$ and $\phi=\varphi_t$ to obtain the estimate for $\dm\n \varphi\dm_2$, we arrive at the relations:
$$\mbox{$\frac12$}\dm \varphi\dm_2^2+(1-\eta)\!\!\intll0t\!\dm\n \varphi\dm_2^2d\tau \leq \mbox{$\frac12$}\dm \varphi_\circ\dm_2^2+  c\!\!\intll0t\! \Big[\dm V\dm_\infty^2+\!  \mbox{${\underset{h=1}{\overset\nu\sum}}$}\dm\n W^h\dm_{r_h}^{\frac{r_h}{r_h-1}}\Big]\dm \varphi\dm_2^2 d\tau
\,,$$
$$\ba{l}\displ\mbox{$\frac12$}\dm\n \varphi\dm_2^2+(1-\eta)\!\!\intll0t\dm P\Delta \varphi\dm_2^2d\tau \leq \mbox{$\frac12$}\dm \n\varphi_\circ\dm_2^2+c\!\!\intll0t\Big[\dm V\dm_\infty^2\dm\n\varphi\dm_2^2+\\\displ\hskip6.5cm \mbox{${\underset{h=1}{\overset\nu\sum}}$}\dm\n W^h\dm_{r_h}^2 \dm\n\varphi\dm_2^{2a_h}\dm\varphi\dm_2^{{2(1-a_h)}}\Big] d\tau 
\,,\ea$$
$$\ba{l}\displ\mbox{$\frac12$}\dm\n \varphi\dm_2^2+(1-\eta)\!\!\intll0t\dm  \varphi_t\dm_2^2d\tau  \leq \mbox{$\frac12$}\dm \n\varphi_\circ\dm_2^2 +c\!\!\intll0t\Big[\dm V\dm_\infty^2\dm\n\varphi\dm_2^2\\\displ\hskip6.3cm+ \mbox{${\underset{h=1}{\overset\nu\sum}}$}\dm\n W^h\dm_{r_h} ^2\dm\n\varphi\dm_2^{2a_h}\dm\varphi\dm_2^{2(1-a_h)}\Big] d\tau
\,.\ea$$ It is known, these estimates allow us to obtain a solution to problem \rf{LNSQ-II} in the $L^2$-setting.  Thanks to the properties of mollifier and the previous result of existence in $L^2$-setting,   for all $q\in(1,2)$, we can consider, not uniformly in $\vep>0$, the right hand side of \rf{LNSQ-II} belonging to $L^q(0,T;L^q(\R^2_+))$.
Hence by virtue of Theorem\,\ref{LPSP}, for all $\vep>0$ we obtain a solution $(\varphi,\pi_\varphi)$ verifying also estimate \rf{LEPSP} with the right hand side of  equation \rf{LNSQ-II}$_1$ in place of $G$, that is for all $s>0$ we get
\be\label{L}\ba{l}\displ\intll0s\Big[\dm \varphi_t(t)\dm_q^q+\dm D^2\varphi(t)\dm_q^q+\dm \n\pi_\varphi(t)\dm_q^q\Big]dt\\\displ\hskip1.4cm\leq c\Big[\dm \varphi_\circ\dm_{{2-\frac2q,q}}^q\!+\!\intll0s\!\dm V\cdot J_\vep[\chi_\vep\n \varphi]+J_\vep[\varphi]\cdot(\n\mbox{${\underset{h=1}{\overset\nu\sum}}$}W^h)^T \dm_q^qdt\Big].\ea\ee Now we look for estimates of $(\varphi,\pi_\varphi)$ in the $L^q$-setting, $q\in(1,2)$, and uniformly with respect to $\vep$.
Applying estimates \rf{UESP-I} to the right hand side of \rf{L}, we get
$$\ba{l}\displ\intll0s\!\Big[\dm \varphi_t(t)\dm_q^q+\dm D^2\varphi(t)\dm_q^q+\dm \n\pi_\varphi(t)\dm_q^q\Big]\!dt\\\displ\hskip3cm\leq c\dm \varphi_\circ\dm_{{2-\frac2q,q}}^q+c\intll0s
\dm V\dm_\infty^q\dm D^2\varphi\dm_q^\frac q2\Big[\intll0t\dm\varphi_\tau\dm_qd\tau+\dm \varphi_\circ\dm_q\Big]^\frac q2dt\\\hskip4.5cm\displ +c \mbox{${\underset{h=1}{\overset\nu\sum}}$} \intll0s \dm\n W^h\dm_{r_h}^q\dm D^2\varphi\dm_q^\frac q{r_h\!\!\null}\Big[\intll0t\dm \varphi_\tau\dm_qd\tau+\dm \varphi_\circ\dm_q\Big]^{\frac q{r_h'\!\!\null}}dt\,.\ea$$
Employing the Cauchy inequality, by means of H\"older's inequality, we deduce
$$\ba{l}\displ\intll0s\!\Big[\dm \varphi_t(t)\dm_q^q+\dm D^2\varphi(t)\dm_q^q+\dm \n\pi_\varphi(t)\dm_q^q\Big]\!dt\\\displ\hskip1.5cm\leq c\dm \varphi_\circ\dm_{{2-\frac2q,q}}^q\!+c\sup_{(0,T)}\Big[\dm V\dm_\infty^{2q}   
 +\!\mbox{${\underset{h=1}{\overset\nu\sum}}$}\dm\n W^h\dm_{r_h}^{q{r_h}'}\Big]\intll0st^{q-1}\Big[\intll0t\!\dm \varphi_\tau\dm_q^qd\tau\Big]dt\,.\ea$$ Employing Gronwall's lemma, we get
$$\ba{l}\displ\intll0s\!\Big[\dm \varphi_t(t)\dm_q^q+\dm D^2\varphi(t)\dm_q^q+\dm \n\pi_\varphi(t)\dm_q^q\Big]\!dt\\\displ\hskip4.5cm\leq c\dm \varphi_\circ\dm_{{2-\frac2q,q}}^q \exp\Big[c\sup_{(0,T)}\Big[\dm V\dm_\infty^{2q}   
 +\mbox{${\underset{h=1}{\overset\nu\sum}}$}\dm\n W^h\dm_{r_h}^{qr_h'}\Big]s^{q}\Big]\ea $$
Since the last inequality is uniform with respect to $\vep>0$, we have  proved that for all $q\in(1,2)$ there exists a solution $(\varphi,\pi_\varphi)$ to problem  \rf{LNSQ-II}. Now we prove that, for any pair $q,\,\ov q\in(1,2)$, the corresponding solutions to problem with the same initial data $\varphi_\circ\in\mathscr C_0(\R^2_+)$ coincide. For this goal we denote by $(\psi,\pi_\psi)$ the difference of the solutions $(\varphi_q,\pi_{\varphi_q})$ and $(\varphi_{\ov q},\pi_{\varphi_{\ov q}})$ corresponding to the same initial data $\varphi_\circ\in \mathscr C_0(\R^2_+)$. Since problem \rf{LNSQ} is linear, $(\psi,\pi_\psi)$ satisfies the same  \rf{LNSQ}, but with all the homogeneous data. The uniqueness is achieved employing the so called weighted function method in this connection see. e.g. \cite{GR}. Multiplying by $\psi e^{-\mu|x|},\mu>0$,    the first equation of \rf{LNSQ} related to $\psi$, setting $g:=e^{-\frac\mu2 |x|}$\,,  we get
$$\ba{l}\frac12\frac d{dt}\dm g\psi \dm_2^2+\dm g\n \psi\dm_2^2=-(\n g^2,\n \psi\cdot\psi)+(V\cdot\n\psi,g^2\psi)\VS\hskip2,3cm+(\psi\cdot\n \mbox{${\overset\nu{\underset{h=1}\sum}}$}W^h,g^2\psi)+(\pi_{\varphi_q}\n g^2,\psi)-(\pi_{\varphi_{\ov q}}\n g^2,\psi)=\mbox{${\overset5{\underset{i=1}\sum}}$}J_i(t)\,.\ea$$
Applying the Cauchy inequality, we get
$$|J_1(t)+J_2(t)|\leq c(\mu+\sup_{(0,T)}\dm V\dm_\infty^2)\dm g\psi\dm_2^2+\eta\dm g\n\psi\dm_2^2\,,$$
Applying H\"older's inequality, and subsequently the Cauchy inequality, we get
$$|J_3(t)|\leq \mbox{${\overset\nu{\underset{h=1}\sum}}$}\dm\n W^h\dm_{r_h}\dm g\psi\dm_{2r'_h}^2
$$ Employing the Gagliardo-Nirenberg inequality, and then  the Cauchy inequality, we obtain $$\ba{l}
|J_3(t)|\leq c\mbox{${\overset\nu{\underset{h=1}\sum}}$}\dm\n W^h\dm_{r_h}\dm \n g\psi+g\n\psi\dm_{2}^\frac2{r_h\!\!\!\null}\;\dm g\psi\dm_2^{\frac2{r'_h\!\!\!\null}}\VS\hskip 4cm \leq c\dm g\psi\dm_{2}^2 \mbox{${\overset\nu{\underset{h=1}\sum}}$}\Big[\dm\n W^h\dm_{r_h}+\dm\n W^h\dm_{r_h}^{r'_h}\Big]+\eta\dm g\n \psi\dm_2^2\,. \ea $$Finally, we consider the terms with pressure fields. It is enough to argument on single term,   the discussion for the other term is anologous. Applying H\"older's inequality and Lemma\,\ref{SOB}, we get
$$|J_4(t)|\leq \mu\dm \pi_{\varphi_q}\dm _{\frac{2q}{2-q}}\dm g\dm_{q'}\dm g\psi\dm_2\leq c\mu^{\frac2q-1}\dm \n \pi_{\varphi_q}\dm_q\dm g\psi\dm_2\,.$$ Collecting the estimates for $J_i$, choosing $\eta$ small, we arrive at 
$$\ba{l}\frac d{dt}\dm g\psi\dm_2\leq c\Big[\mu+\sup_{(0,T)}\big(\dm V\dm_\infty^2+\mbox{${\overset\nu{\underset{h=1}\sum}}$}\big[\dm\n W^h\dm_{r_h}+\dm\n W^h\dm_{r_h}^{r'_h}\big]\big)\Big]\dm g\psi\dm_2\VS\hskip7cm+c(\mu^{\frac2{q }-1}\dm \n\pi_{\varphi_q}\dm_q+\mu^{\frac2{\ov q }-1}\dm\n \pi_{\varphi_{\ov q}}\dm_{\ov q})\,.\ea$$ Integrating the last differential inequality,   uniformly in $\mu>0$, we deduce
$$\dm g\psi\dm_2\leq c\exp{t\, \frak C(t)}\intll0t (\mu^{\frac2{q}-1}\dm \n\pi_{\varphi_q}\dm_q+\mu^{\frac2{\ov q }-1}\dm\n \pi_{\varphi_{\ov q}}\dm_{\ov q})dt\,,$$ where we set $\frak C(t):=c\Big[\mu+\sup_{(0,T)}\big(\dm V\dm_\infty^2+\mbox{${\overset\nu{\underset{h=1}\sum}}$}\big[\dm\n W^h\dm_{r_h}+\dm\n W^h\dm_{r_h}^{r'_h}\big]\big)\Big]$. Hence in the limit for $\mu\to0$ we prove the uniqueness. 
 \ep
 
\section{A special IBVP of the perturbed Navier-Stokes equations}
In this section we study the following initial boundary value problem:
\be\label{PNS}\ba{l}w_t-\Delta w-\!\n \pi_{w}=-w\cdot\n w-V\cdot\!\n w-w\cdot\!\n V+F \,,\VS\n\!\cdot\! w=0\mbox{ in }(0,T)\!\times\!\R^2_+,\VS w=0\mbox{ on }(0,T)\times\{x_2=0\},\quad w=0\mbox{ on }\{0\}\times\R^2_+.\ea\ee
The following result holds:
\begin{tho}\label{TPNS}{\sl Assume that $V\in L^\infty((0,T)\times \R^2_+)$, satisfying the divergence free, and, for some $r\geq2$, $\n V\in L^\infty(0,T;L^r(\R^2_+))\cap L^r(0,T;W^{1,r}(\R^2_+))$. If $F\in L^2(0,T;L^2(\R^2_+)\cap L^q(0,T; L^{q}(\R^2_+))$, for some $q>2$,   then there exists a unique solution to problem \rf{PNS} such that for all $T>0$\be\label{TPNS-I}\ba{c}w\!\in\! C([0,T);J^{1,2}(\R^2_+)) \,,\, D^2_{x_ix_j}w\,,\,w_t\,,\!\n\pi_w\!\in\! L^2(0,T;L^2(\R^2_+)), \VS w\in C([0,T)\times\ov\R^2_+)\mbox{\, and\, }\lim_{t\to0}\dm w(t) \dm_\infty=0\,,\VS\pi_w(t,x)-\pi_w(t)=j(t)O(|x|^\gamma\!+1)\mbox{\;a.e.\,in }t>0,\mbox{\,for all\;}x\!\in\!\R^2_+\,,  \ea\ee where we have set $\pi_w(t):=\pi_w(t,0)$,   $j(t)\in L^{{\ov\sigma}}(0,T)$,     $\gamma\in(0,1)$ and $\ov \sigma  >2$.}\end{tho}
\bp We  set $\frak D(t):= {\rm exp}\Big[t+c \!\intll0t\! \dm\n V\dm_r^{r'}\,  d\tau\Big] \intll0t\dm F\dm_2^2\;d\tau$. By making use of the Galerkin method, essentially in the way employed in  section 4 for problem \rf{LNS}, we can deduce the following inequalities:
\be\label{PEI}\ba{l}\displ\dm w(t)\dm_2^2+\intll0t\dm\n w\dm_2^2d\tau\leq  \frak D(t)\,,\ea\ee
\be\label{PEII}\ba{l}\displ\dm \n  w(t)\dm_2^2+\intll0t\Big[\dm w_\tau\dm_2^2+\dm P\Delta w\dm_2^2\Big]d\tau\VS\hskip1.8cm\leq c\exp[\frak D^2(t)]\Big[ t\frak D(t)\sup_{(0,T)}\big(\dm V\dm_\infty^4\!+\dm \n V\dm_{r}^{2r'}\big)\!+ \!\! \intll0t\!\dm F\dm_2^2d\tau\Big], \ea\ee 
uniformly in $t>0$. Hence we consider as achieved \rf{TPNS-I}$_{1}$\,. We look for the following decomposition of the solution: $$w:=w^1+w^2+w^3\mbox{\; and \; }\pi_w:=\pi_{w^1}+\pi_{w^2}+\pi_{w^3}\,,$$ where $(w^1,\pi_{w_1})$ is the solution to problem \rf{ST} with zero initial data and force data $G^1:=-(V+w)\cdot\n w$, further $(w^2,\pi_{w_2})$ is the solution to problem \rf{ST} with zero initial data and $G^2:=-w\cdot\n V$\,, and finally $(w^3,\pi_{w^3})$ is the solution to problem \rf{ST} with zero initial data and $G^3:=F$.  
By virtue of estimate \rf{PEI}-\rf{PEII}, employing the Gagliardo-Nirenberg inequality, we get $G^1\in L^\frac83(0,T;L^\frac83(\R^2_+))$. Hence Theorem\,\ref{LPSP} ensures the existence of a unique solution such that $w^1\in C(0,T;J^{1,\frac83}(\R^2_+))\cap L^\frac83(0,T;W^{2,\frac83}(\R^2_+))$, and by Sobolev embedding  theorem $w^1\in C([0,T);C(\ov \R^2_+))$. We also get $\n\pi_{w^1}\in L^\frac83(0,T;L^\frac83(\R^2_+))$. Analogously, by virtue of \rf{PEI}-\rf{PEII} under our assumption for $V$, for some $r_1>2$,  we get $G^2\in L^{ r_1}(0,T;L^{r_1}(\R^2_+))$.   Hence Theorem\,\ref{LPSP} ensures the existence of a unique solution such that $w^2\in C(0,T;J^{1,r_1}(\R^2_+))\cap L^{r_1}(0,T;W^{2,r_1}(\R^2_+))$, and by Sobolev embedding  theorem $w^2\in C([0,T);C(\ov \R^2_+))$. We also get $\n\pi_{w^2}\in L^{r_1}(0,T;L^{r_1}(\R^2_+))$. Finally, by hypotheses on $F$ we have $G^3\in L^q(0,T;L^q(\R^2_+))$ for some $q>2$. 
Hence Theorem\,\ref{LPSP} ensures the existence of a unique solution such that $w^3\in C(0,T;J^{1,q}(\R^2_+))\cap L^{q}(0,T;W^{2,q}(\R^2_+))$, and by Sobolev embedding  theorem $w^3\in C([0,T);C(\ov \R^2_+))$. We also get $\n\pi_{w^3}\in L^{q}(0,T;L^{q}(\R^2_+))$. As well as for  the solutions we get \be\label{ICI-I}\lim_{t\to0}\dm w^1(t)\dm_\infty=\lim_{t\to0}\dm w^2(t)\dm_\infty=\lim_{t\to0}\dm w^3(t)\dm_\infty=0\,.\ee The difference $(w-w^1-w^2-w^3, \pi_w-\pi_{w^1}-\pi_{w^2}-\pi_{w^3})$ is a solution to the Stokes problem with homogenous data. Hence it is easy to prove that, up to a constant for the pressure field, the difference is identically zero. Therefore the decomposition holds and we  deduce \rf{TPNS-I}$_{2}$ via \rf{ICI-I} and \rf{TPNS-I}$_3$ employing Sobolev embedding theorem and setting $j(t):=\dm \n \pi_{w^1}(t)\dm_\frac83+\dm\n \pi_{w^2}(t)\dm_{r_1}+\dm \n\pi_{w^3}(t)\dm_q$. Finally, setting $\gamma:=\min\{\frac14,1-\frac2{r_1},1-\frac2q\}$, as well setting $\ov \sigma=\min\{\frac83,r_1,q\}$, we complete the proof.\ep
\section{Proof of Theorem\,\ref{CT}}\subsection{\label{EX}Existence} We develop the proof distinguishing the following cases for the initial data:
\begin{itemize}\item[1)] $p=2$\,,\item[2)]  $p\in(2,4]$\,,\item[3)]  $p>4$\,.\end{itemize} 
In all the cases 1)-2) by $(U,\pi_U)$ we denote the solution to problem \rf{ST} assuming data $U(0,x):=u_0$ whose existence is ensured by Theorem\,\ref{GPST}.\vskip0.15cm\par 1) 
In the case $p=2$,   we consider problem \rf{PNS} with coefficient $V:=U$ and data force $F=U\cdot\n U$. This data $F\in L^\infty(0,T;L^2(\R^2_+))\cap L^\frac83(0,T;L^\frac83(\OO))$\,\footnote{\,For the last claim we apply the interpolation of Lebesgue spaces, hence we get$$\dm \n U(t)\dm_\frac83\leq c\dm \n U(t)\dm_\infty^\frac14\dm \n U(t)\dm_2^\frac34\,,\mbox{ for all }t>0\,.$$ This last estimate and  \rf{ESP-I} imply $U\cdot\n U\in L^\frac83(0,T;L^\frac83(\OO))$.}. Denoted by $(w,\pi_w)$ the solution of Theorem\,\ref{TPNS}, setting $u:=U+w$ and $\pi_u:=\pi_U+\pi_w$ we have proved the existence.\vskip0.15cm\par 2) In  the case of $p\in(2,4]$, we denote by $(v^1,\pi_{v^1})$ the solution to problem \rf{LNS} where we assume  for data:
$$h=0, \;w^0:=U\quad\mbox{and }\quad F:=U\cdot \n U\,.$$  By the regularity of $U$, for all $T>0$, we have that $U\cdot\n U\in L^p(0,T;L^p(\R^2_+))$. By virtue of Theorem\,\ref{TLNS}, for all $T>0$, we have $v^1\in C([0,T);J^{1,p}(\R^2_+))$. Since $p\in (2,4]$, we obtain $ v_1\in C(0,T;L^q(\R^2_+))$ for all $q\in[p,\infty]$. Hence we also get that \be\label{LII}\mbox{for all }T>0,\quad v^1\cdot\n v^1\in L^2(0,T;L^2(\R^2_+))\cap L^{\frac83}(0,T;L^\frac83(\R^2_+))\,.\ee 
Actually, if $p\in(2,4]$, applying H\"older's inequality,   we get $$\ba{l}\dm v^1\cdot\n v^1\dm_2\leq\dm v^1\dm_{\frac{2p}{p-2}}\dm\n v^1\dm_p<\infty,\mbox{ for all }t>0   \,,\VS \dm v^1\cdot\n v^1\dm_p\leq \dm v^1\dm_\infty\dm \n v^1\dm_p<\infty,\mbox{ for all }t>0\, .\ea$$  The above estimates justify \rf{LII}. Now we introduce the following perturbed Navier-Stokes problem:
\be\label{APNS}\ba{l}v^2_t-\Delta v^2+\!\n \pi_{v^2}=-v^2\cdot\n v^2-\mbox{${\overset{1}{\underset{\ell=0}\sum}}$}w^\ell\!\cdot\!\n v^2-v^2\cdot\!\n \mbox{${\overset{1}{\underset{\ell=0}\sum}}$}w^\ell\!+F^1\,,\VS\n\cdot v^2=0\mbox{ in }(0,T)\times\R^2_+,\VS v^2=0\mbox{ on }(0,T)\times\{x_2=0\},\quad v^2=0\mbox{ on }\{0\}\times\R^2_+,\ea\ee
where $w^0:=U,\;w^1:=v^1$ and $F^1:=-v^1\cdot\n v^1$. By virtue of Theorem\,\ref{TPNS}, there exists a unique regular solution $(v^2,\pi_{v^2})$ to problem \rf{APNS}. Hence setting $u:=U+v^1+v^2$ and $\pi_u:=\pi_U+\pi_{v^1}+\pi_{v^2}$ we have proved the existence claimed in   Theorem\,\ref{CT} for  $p\in(2,4]$.\vskip0.15cm \par 3)
In the case of $p>4$,  we consider the greatest integer floor $  k$ of $\log_2\frac p2$   such that $\frac p{2^{  k}\!\!\!\null}\,> 2$ and
$\frac p{2^{  k+1}\hskip-0.4cm\null}\hskip0.3cm\leq2$. Hence we get $\frac p{2^{  k}\!\!\!\null}\hskip0.1cm\in(2,4]$. 
We consider the finite family of solutions $\{ (v^h,\pi_{v^h})\}_{h=1,\dots,k+1}$ to problem \rf{LNS}, where
$$\ba{l}\mbox{we set }w^0:=U\,\mbox{ with initial data }u_0\,, \mbox{and,  for all }\ell=1,\dots,k+1,\VS X^\ell:=C(0,T;J^{1,\frac p{2^{\ell-1\!\!\!\!\!\!\!\!\null}}}\hskip0.3cm(\R^2_+))\cap L^\frac p{2^{\ell-1\!\!\!\!\!\!\!\!\null}}\hskip0.3cm(0,T; W^{2,\frac p{2^{\ell-1\!\!\!\!\!\!\!\!\null}}}\hskip0.3cm(\R^2_+))\,,\VS v^1\,\mbox{ has as coefficient } w^0=U \mbox{ and as force data}\VS \hskip 4cmF^1:=-U\cdot\n U\in L^p(0,T;L^p(\R^2_+))\,,\VS
v^2\,\mbox{ has as coefficients }w^0=U\,,\;w^1=v^1\,,\mbox{ with }w^1\in  X^1\,, \mbox{ as force data}\VS\hskip4cmF^2:=-v^1\cdot\n v^1\in L^\frac p2(0,T;L^\frac p2(\R^2_+))\,,\VS\cdots\quad\cdots\quad\cdots\VS v^{k+1}\mbox{ has as coefficients }w^0=U\,,w^1=v^1\,,\cdots,w^{k}=v^{k}, \;w^\ell\in  X^\ell\,,\,\ell =1,\dots,k\,, \VS\mbox{as force data} \VS \hskip4cmF^{k+1}:=-v^{k}\cdot\n v^{k}\in L^{\frac p{2^{k \!\!\!\null}}}\hskip0.1cm(0,T;L^{\frac p {2^{k  \!\!\!\null}}}\hskip0.1cm(\R^2_+))\,.\ea $$
For each $h=1,\dots, k+1$ the existence of the pair $(v^h,\pi_{v^h})$ is ensured by Theorem\,\ref{TLNS}. Now, by means of Theorem\,\ref{TPNS} we can solve the problem
\be\label{KPNS}\ba{l}w_t-\Delta w+\!\n \pi_{w}=-w\cdot\n w-\mbox{${\overset{k+1}{\underset{\ell=0}\sum}}$}w^\ell\!\cdot\!\n w-w\cdot\!\n \mbox{${\overset{k+1}{\underset{\ell=0}\sum}}$}w^\ell\!+F \,,\VS\n\!\cdot\! w=0\mbox{ in }(0,T)\!\times\!\R^2_+,\VS w=0\mbox{ on }(0,T)\times\{x_2=0\},\quad w=0\mbox{ on }\{0\}\times\R^2_+,\ea\ee
where $w^0\!:=U,\, w^\ell\!:=v^\ell,\ell=1,\dots,k+1\,,$ and $F\!:=-v^{k+1}\cdot\n v^{k+1}\!\in\! L^2(0,T;L^2(\R^2_+))$.
The last claim is consequence of the fact that $\frac p{2^k\!\!\!\null}\in(2,4]$\,.
Setting $u:=U+w+\mbox{${\overset{k+1}{\underset{\ell=1}\sum}}$}v^\ell$ and $\pi_u:=\pi_U+\pi_w+\mbox{${\overset{k+1}{\underset{\ell=0}\sum}}$}\pi_{v^\ell}$, then by construction the pair $(u,\pi_u)$ is  a solution to equations \rf{NS}. The solution satisfies the boundary condition.   Since $\lim_{t\to0}U(t,x)=u_0(x)$ for all $x\in\R^2_+$, in order to prove that $\lim_{t\to0}u(t,x)=u_0(x)$ for all $x\in \R^2_+$,
we can limit ourselves to prove that $\lim_{t\to0}\dm w+\mbox{${\overset{k+1}{\underset{\ell=1}\sum}}$}v^\ell\dm_\infty=0$. Actually, this is a consequence of the estimates \rf{TLNS-I}$_1$ for the functions $v^\ell$, $\ell=1,\dots,k+1$, and 
of estimate \rf {TPNS-I}$_2$ for $w$.  The proof of the existence is completed if the regularity properties 
hold. This is a classical result for solutions to the 2D-Navier-Stokes  problem, hence it is omitted.
\subsection{Uniqueness}We begin recalling that the pressure field $\pi_u$ determined in section\,\ref{EX} is given by the sum $$\pi_u:=\pi_U+\mbox{$\underset{\ell=1}{\overset {k+1}\sum}$}\pi_{v^\ell}+\pi_w\,.$$ For each term we have the following estimates:\begin{itemize}\item[iv.] for all $\gamma\in(0,1)$, $|\pi_U(t,x)|=O(|x|^{\gamma})t^{-\frac\gamma2}$ for all $(t,x)\in (0,T)\times\R^2_+$\,,
\item[v.] for all $\ell=1,\dots,k+1$, almost everywhere in $t>0$, $\nabla\pi_{v^\ell}\in L^{\frac p {2^{\ell \!\!\!\null}}}\,(\R^2_+)$ that, via Sobolev embedding,   furnishes    $$|\pi_{v^\ell}(t,x)- \pi_{v^\ell}^0(t)|\leq c J_\ell(t) |x-x^0|^{1-\frac {2^\ell\hskip-0.1cm\null} p}\hskip0.2cm,  $$ almost everywhere in $t>0$ and for all $x\in \R^2_+$\,, where we set $\pi_{v^\ell}^0(t):=\pi_{v^\ell}(t,x_0)$ and  $J_\ell(t):=\dm\n \pi_{v^\ell}(t)\dm_{\frac p{2^\ell\hskip-0.1cm\null}}\,,$ since $ 1-\frac{2^{k+1}\hskip-0.4cm\null}p\hskip0.3cm\leq 1-\frac{2^{\ell\hskip-0.1cm\null} }p<1$, for all $\ell=1,\dots, k+1$, then,  for a suitable constant $c$, we get   $$|\pi_{v^\ell}(t,x)-\pi^0_{v^\ell}(t)| \leq cJ(t)\big[|x |^{1-\frac{2^{k+1\!\!\! \!\!\!\!\!\!\null}}p}\hskip0.2cm+1\big]\,,$$ almost everywhere in $t>0$ and for all $x\in R^2_+$, where we have set $J(t):={\underset{\ell=1}{\overset{k+1}\sum}}\, \dm \n \pi_{v^\ell}(t)\dm_{\frac p {2^{\ell \!\!\!\null}}}$ which belongs to $  L^{\frac p{2^{k+1\!\!\!\!\!\!\!\!\null}}}\hskip0.3cm(0,T) \,,$\item[vi.] finally, for $\pi_w(t,x)$ estimate \rf{TPNS-I}$_3$ holds.  \end{itemize}The following lemma concerns  a weighted energy estimate of the kind proved in \cite{GR}.
    \begin{lemma}\label{LU}{\sl Let $(u,\pi_u)$ and $(\ov u,\pi_{\ov u})$ be two regular solutions to problem \rf{NS} assuming the same initial data $u_0\in L^\infty(\R^2_+)$. Assume that the fields $u,\,\ov u\in L^\infty((0,T)\times\R^2_+)$. Assume that there exist $\rho,\,\ov \rho\in (0,1)$  such that\be\label{LU-I}\ba{l}|\pi_u(t,x)-\pi_u(t)|\leq J_u(t)(|x|+1)^\rho,\VS|\pi_{\ov u}(t,x)-\pi_{\ov u}(t)|\leq J_{\ov u}(t)(|x|+1)^{\ov\rho},\;t>0,\;x\in\R^2_+,\ea\ee where functions  $\pi_u(t),\pi_{\ov u}(t),J_u(t),J_{\ov u}(t)\in L^1_{\ell oc}[0,T)$. Then, for $\beta\in (\rho,1)\cap(\ov\rho,1)$, we get
\be\label{LU-II}\dm (u(t)-\ov u(t))(|x|+1)^\beta\dm_2<\infty,\mbox{ for all }t\in (0,T)\,.\ee
}\end{lemma}\bp We set $W:=u-\ov u$ and $\pi_W:=\pi_u-\pi_{\ov u}+\pi_{\ov u}(t)-\pi_u(t)$. Denoted by $\rho_0:=\max\{\rho,\ov \rho\}$, by the assumption \rf{LU-I} we have  $|\pi_W(t,x) | \leq (J_u(t)+J_{\ov u}(t))(|x|+1)^{\rho_0},t>0,x\in \R^2_+$. We consider the system of the difference $(W,\pi_W)$ written as:
$$\ba{l}W_t-\Delta W+\n\pi_W=-W\cdot\n u-\ov u\cdot \n W\,,\quad \n\cdot W=0\mbox{ in }(0,T)\times\R^2_+\VS W=0\mbox{ on }(0,T)\times\{x_2=0\},\quad W=0\mbox{ on }\{0\}\times\R^2_+ \,.\ea$$
We multiply the first equation by $(|x|+1)^{-2\beta}\exp[-2\mu |x|]W=:g^2(x)W$. Integrating by parts on $(0,t)\times\R^2_+$
 we get the weighted energy relation $$ \mbox{\large$\frac12$} \dm W(t)g\dm_2^2+\intll0t\dm g\n W(\tau)\dm_2^2d\tau= \intll0t\mbox{${\overset4{\underset{i=1}\sum}}$}I_i(t)dt\,.$$
Integrating by parts and applying H\"older's inequality, we get
\begin{itemize}\item$\displ|I_1(t)|= |\intl{R^2_+}\n g^2\cdot\n W\cdot W dx\leq 2(\beta+\mu)\dm Wg\dm_2\dm g\n W\dm_2
\,,$\item$\displ|I_2(t)|=|\intl{\R^2_+}W\cdot\n W\cdot ug^2dx+\intl{\R^2_+}u\cdot WW\cdot\n g^2 dx|$\newline $\null\hskip4cm\leq \dm u\dm_\infty\Big[\dm Wg\dm_2\dm \n Wg\dm_2+2(\beta+\mu)\dm Wg\dm_2^2\Big]\,,$\item$\displ |I_3(t)|=\frac12|\intl{\R^2_+}|W|^2\ov u\cdot\n g^2 dx| \leq (\beta+\mu)\dm \ov u\dm_\infty\dm Wg\dm_2^2\,,$ \item $\displ|I_4(t)|=|\intl{\R^2_+}\pi_W\n g^2\cdot Wdx|$\newline$\null\hskip2.3cm\leq 2(J_u(t)+J_{\ov u}(t))\Big[\beta\dm\frac{(|x|+1)^{\rho_0}\!\!\!\null}{(|x|+1)}\hskip0.1cmg\dm_2+\mu\dm(1+|x|)^{\rho_0} g\dm_2\Big]\dm Wg\dm_2\,.$\end{itemize} We collect the estimates for $I_i$, and we estimate the right hand side of the energy relation. After applying the Cauchy inequality,  we get $$\ba{l}\displ\dm W(t)g\dm_2\leq \intll0t\Big[C(\beta,\mu)+\dm u\dm_\infty^2+C(\beta,\mu)\big[\dm u\dm_\infty+\dm\ov u\dm_\infty\big]\Big]\dm Wg\dm_2d\tau\\\hskip2.3cm\displ+2\Big[\beta\dm\frac{(|x|+1)^{\rho_0}\!\!\!\null}{(|x|+1)}\hskip0.1cmg\dm_2+\mu\dm(1+|x|)^{\rho_0} g\dm_2\Big]\intll0t (J_u(\tau)+J_{\ov u}(\tau))d\tau\,.\ea$$
Employing Gronwall lemma we get
\be\label{LU-III}\dm W(t)g\dm_2\leq M(\beta,\mu) H(t,\beta,\mu)\!\!\intll0t\!(J_u(\tau)+J_{\ov u}(\tau))d\tau\,,\ee where we have set $H(t,\beta,\mu):=e^{\intll0t\big[C(\beta,\mu)+\dm u\dm_\infty^2+c(\beta,\mu)(\dm u\dm_\infty+\dm\ov u\dm_\infty)\big]d\tau}$, and $$M(\beta,\mu):=2\Big[\beta\dm\frac{(|x|+1)^{\rho_0}\!\!\!\null}{(|x|+1)}\hskip0.1cmg\dm_2+\mu\dm(1+|x|)^{\rho_0} g\dm_2\Big]\,.$$ Since $\beta>\rho_0$,   a simple computation furnishes that
$$M(\beta,\mu)\leq c\dm{(|x|+1)^{\rho_0-1-\beta} } \dm_2=:\ov M<\infty\,,\quad\mbox{for all }\mu>0.$$ Hence, \rf{LU-III} becomes
$$\dm W(t)g\dm_2\leq \ov MH(t,\beta,\mu)\intll0t(J_u(\tau)+J_{\ov u}(\tau))d\tau\,.$$ Since $H(t,\beta,\mu)\leq H(t,\beta,1)$, letting $\mu\to0$,  via the Beppo Levi theorem,   the last estimate leads to the thesis.\ep
Now we are in a position to obtain the uniqueness. We employ a duality argument, which is  a variant of the one introduced in \cite{F} for the uniqueness to solution of the Navier-Stokes Cauchy problem.  Actually, the result is a consequence of the following:
\begin{lemma}{\sl Let $(u,\pi_u)$ be the solution to problem \rf{NS} furnished by section\,\ref{EX}. Let $(\ov u,\pi_{\ov u})$ be a regular solution to problem \rf{NS} corresponding to the same data $u_0$. Assume that $\ov u\in L^\infty((0,T)\times\R^2_+)$. Further assume that $\pi_{\ov u}$ satisfies condition \rf{LU-I}$_2$. Then, up to a function $c(t)$, the given solutions coincide on $(0,T)$.}\end{lemma} \bp We denote by $\widehat V:=u-\ov u$ and $\widehat\pi:=\pi_u-\pi_{\ov u}$. The pair $(\widehat V,\widehat\pi)$ is a solution to the problem
\be\label{V}\ba{l}\widehat V_t-\Delta \widehat V+\n\widehat\pi=-\ov u\cdot\n \widehat V-\widehat V\cdot \n u\,,\quad \n\cdot \widehat V=0\mbox{ in }(0,T)\times\R^2_+\,,\VS \widehat V=0\mbox{ on }(0,T)\times\{x_2=0\},\quad \widehat V=0\mbox{ on }\{0\}\times\R^2_+\,.\ea\ee
We start claiming that $\widehat V $ enjoys of estimate \rf{LU-II}. For this it is enough to verify the hypotheses of Lemma\,\ref{LU}. Of course, we can limit ourselves to verify assumption \rf{LU-I}$_1$. This assumption is a consequence of items iv.-vi.\,. Actually, a simple computation allows us to say that $\rho=1-\frac {2^{k+1}\!\!\!\!\!\!\!\!\null}p$\hskip0.33cm in \rf{LU-I}$_1$. We also remark that by construction $\n u={\underset{\ell=1}{\overset{k+1}\sum}}\n v^\ell+\n U+\n w=:{\underset{h=1}{\overset{k+3}\sum}}\n u^h\,,$ hence each term of the sum belongs to $C(0,T;L^{r_h}(\R^2_+))$ for a suitable $r_h>2$. Thanks to Theorem\,\ref{TLNSQ}, for a given function $\varphi_\circ\in \mathscr C_0(\R^2_+)$, we can state the existence of the  solutions $(\varphi,\pi_{\varphi})$ to problem \rf{LNSQ} with $-\widehat{\ov u}$ in place of $V$ and $\mbox{${\underset{h=1}{\overset{k+3}\sum}}$}\widehat u^h$ in place of $W$, where for a fixed $t>0$ we have set $$\widehat{\ov u}(\tau,x):=\ov u(t-\tau,x)\mbox{ and }\widehat u^h(\tau,x):=u^h(t-\tau,x),\mbox{ for all }\tau\in[0,t]\,.$$ In such a way the pair $$\widehat\varphi(\tau,x):=\varphi(t-\tau,x)\quad\widehat\pi_{\varphi}(\tau,x):=\pi_\varphi(t-\tau,x)\,,\mbox{ for all }\tau\in[0,t]\,,$$is a solution to the adjoint problem of \rf{V}:
\be\label{VA}\ba{l}\widehat\varphi_\tau\!+\!\Delta\widehat\varphi+\!\hskip-0.05cm\n\pi_{\widehat\varphi}=\!-\ov u\cdot\!\n\widehat\varphi+\!\widehat\varphi\cdot\!({\underset{h=1}{\overset{k+3}\sum}}\hskip-0.05cm\n u^h)^T,\;\;\n\!\cdot\!\widehat\varphi=0,\mbox{\,in\,}(0,t)\!\times\!\R^2_+ ,\VS  \widehat\varphi=0\mbox{ on }(0,T)\times\{x_2=0\},\quad \widehat\varphi=\varphi_\circ\mbox{ on }\{t\}\times\R^2_+.\ea\ee
We multiply equation \rf{V}$_1$ by $\chi_R(x)\widehat \varphi$. We assume that $R>$diam$\{$supp$\,\varphi_\circ\}$. Integrating by parts on $(0,t)\times\R^2_+$, we get
$$\ba{l}\displ(\widehat V\!,\varphi_\circ)=\!\!
\intll0t\!\Big[(\widehat V\chi_R,\widehat \varphi_\tau\!+\Delta\widehat \varphi)\!+\!(\widehat V\!,\Delta\chi_R\widehat \varphi)\!+\!2(\widehat V\!,\n\chi_R\!\cdot\!\n\widehat\varphi)
\!+\! (\ov u\cdot\n\widehat \varphi, \widehat V\chi_R)\\\displ\hskip4cm+(\ov u\cdot\!\n\!\chi_R,\widehat V\cdot\widehat \varphi)-(\widehat V\cdot\!\n u,\chi_R\widehat\varphi)  +(\widehat\pi,\n\chi_R\cdot\varphi)\Big]d\tau.\ea$$   In the previous relation substituting the right hand side of \rf{VA}, integrating by parts, we get
\be\label{LE}\ba{ll}\displ(\widehat V,\varphi_\circ)\hskip-0.2cm&\displ=
\intll0t\Big[ (\widehat V,\Delta\chi_R\widehat \varphi)+2(\widehat V,\n\chi_R\cdot\n\widehat\varphi)
 +(\ov u\cdot\n\chi_R,\widehat V\cdot\widehat \varphi)\\&\displ\hskip1.3cm+ (\widehat\pi,\n\chi_R\cdot\varphi)+(\pi_{\widehat\varphi},\n\chi_R\cdot \widehat V)\Big]d\tau=\mbox{${\overset 5{\underset{i=1}\sum}}$}\intll0tJ_i(\tau)d\tau.\ea
\ee Since, $\widehat V,\ov u\in L^{\infty}((0,T)\times\R^2_+)$, and $\widehat\varphi\in C(0,T; L^q(\R^2_+))\cap L^q(0,T;J^{1,q}(\R^2_+))$, for all $q\in(1,2)$, applying H\"older's inequality, we immediately deduce that $$\lim_{R\to\infty}|J_1(t)+J_2(t)+J_3(t)|=0\,,\mbox{\; for all }t\in(0,T).$$
Since $\widehat\pi=\pi_{u}-\pi_{\ov u}$, and $\pi_u$ verifies \rf{LU-I}$_1$ with $\rho=1-\frac{2^{k+1}\!\!\!\!\!\!\!\!\null}p$\;\;\; and $\pi_{\ov u}$ verifies \rf{LU-I}$_2$ by hypothesis, applying H\"older's inequality, we arrive at
$$|J_4(t)|\leq c\hskip-0.3cm\intl{R<|x|<2R}\hskip-0.3cm\Big[|x|^{\rho-1}|\varphi|+|x|^{\ov\rho-1}|\varphi|\Big]dx\leq c(R^{\rho+1-\frac2{r } }\dm \varphi\dm_r+
R^{\ov\rho+1-\frac2{\ov r } }\dm \varphi\dm_{\ov r})\,.$$ Choosing $r,\ov r$ in such a way that the exponent of $R$ are negative, we get that
$$\lim_{R\to\infty}|J_4(t)|=0\,,\mbox{\; for all }t\in (0,T)\,.$$ Analogously, recalling that $\widehat V$ has finite weighted energy \rf{LU-II}, applying H\"older's inequality with exponents  $\frac2\beta\,,\frac2{1-\beta}$ and $2$, we get
$$\ba{ll}|J_5(t)|\hskip-0.2cm&=|(\pi_{\widehat\varphi}-{\pi_{\widehat\varphi}}_0,\n\chi_R\cdot\widehat V)|\VSE\leq \dm \pi_{\widehat\varphi}-{\pi_{\widehat\varphi}}_0\dm_{L^\frac2\beta(R<|x|<2R)}\dm |x|^{\beta-1}\dm_{L^\frac2{1-\beta}(R<|x|<2R)}\dm \widehat V|x|^{-\beta}\dm_{L^2(R<|x|<2R)} \,,\ea$$ here we have introduced the constant ${\pi_{\widehat\varphi}}_0$ of   Lemma\,\ref{SOB}. We remark that the introduction of any constant is allowed by the fact that $(\n\chi_R, \widehat V)=0$. As well we remark that for all $\beta\in(0,1)$ there exists a $r\in(1,2)$ such that $\frac2\beta=\frac{2r}{2-r}$. Concurrently made these remarks   justify the estimate of $J_5$. Thus, employing \rf{SOB-I}, the following holds:
$$|J_5(t)|\leq c \dm \n \pi_{ \widehat\varphi}\dm_r \dm \widehat V|x|^{-\beta}\dm_{L^2(R<|x|<2R)}\,.$$
Hence we deduce
$$\lim_{R\to\infty}|J_5(t)|=0\,,\mbox{\; for all }t\in(0,T)\,.$$ Collecting the estimates related to $J_i(t)$ we have that for $R\to\infty$ the right hand side of \rf{LE} is zero, that proves $$(\widehat V,\varphi_\circ)=0\,,\mbox{\; for all }\varphi_\circ\in \mathscr C_0(\R^2_+)\,.$$ Therefore, function $\widehat V$ is the gradient of some $H(t,x)$, which is harmonic for  $\widehat V$ has divergence free. Hence we can claim that $\widehat V\equiv0$. The uniqueness is proved.
\ep

  \end{document}